\documentclass[final]{siamltex}

\usepackage[dvips]{graphicx}
\usepackage{amssymb,amsfonts,amsmath}
\usepackage{eucal,bm}
\usepackage{gensymb,subfigure}
\usepackage{color}
\usepackage{booktabs}
\usepackage{url}

\newcommand{\cE}{{\cal E}}
\newcommand{\cD}{{\cal D}}
\newcommand{\cG}{{\cal G}}
\newcommand{\cW}{{\cal W}}
\newcommand{\cV}{{\cal V}}

\def\cB{{\cal B}}

\newcommand{\R}{\mathbb{R}}

\newcommand{\bmm}[1]{ {\mbox{\boldmath $ {#1}$}}}

\newcommand{\beq}{\begin{equation}}
\newcommand{\eeq}{\end{equation}}

\newcommand{\beqnr}{\begin{eqnarray}}
\newcommand{\eeqnr}{\end{eqnarray}}

\newcommand{\benum}{\begin{enumerate}}
\newcommand{\eenum}{\end{enumerate}}

\newcommand{\QED}{\rule{.1in}{.1in}}

\newcommand{\cR}{{\cal R}}
\newcommand{\cS}{{\cal S}}
\newcommand{\cM}{{\cal M}}

\newtheorem{DE}{Definition}[section]

\newtheorem{LE}[DE]{Lemma}

\newtheorem{PRO}[DE]{Procedure}

\title{Chance Constrained Optimal Power Flow:\\
Risk-Aware Network Control under Uncertainty}
\author{
Daniel Bienstock
\thanks{Department of Industrial Engineering and Operations Research and
Department of Applied Physics and Applied Mathematics, Columbia University, 500 West 120th St. New York, NY 10027 USA ({\tt dano@columbia.edu}).}
\and
Michael Chertkov
\thanks{Theoretical Division and Center for Nonlinear Studies, Los Alamos National Laboratory, NM 87545 USA ({\tt chertkov@lanl.gov}).}
\and
Sean Harnett
\thanks{Department of Applied Physics and Applied Mathematics, Columbia University, 500 West 120th St. New York, NY 10027 USA and Center for Nonlinear Studies, Los Alamos National Laboratory, NM 87545 USA.}
}

\begin{document}

\maketitle

\begin{abstract}
When uncontrollable resources fluctuate, Optimum Power Flow (OPF), routinely used by the electric power industry to re-dispatch hourly controllable generation (coal, gas and hydro plants) over control areas of transmission networks, can result in grid instability, and, potentially, cascading outages. This risk arises because OPF dispatch is computed without awareness of major uncertainty, in particular fluctuations in renewable output.  As a result, grid operation under OPF with renewable variability can lead to frequent conditions where power line flow ratings are significantly exceeded. Such a condition, which is borne by simulations of real grids, would likely resulting in automatic line tripping to protect lines from thermal stress, a risky and undesirable outcome which compromises stability. Smart grid goals include a commitment to large penetration of highly fluctuating renewables, thus calling to reconsider current practices, in particular the use of standard OPF. Our Chance Constrained (CC) OPF corrects the problem and mitigates dangerous renewable fluctuations with minimal changes in the current operational procedure. Assuming availability of a reliable wind forecast parameterizing the distribution function of the uncertain generation, our CC-OPF satisfies all the constraints with high probability while simultaneously minimizing the cost of economic re-dispatch. CC-OPF allows efficient implementation, e.g. solving a typical instance over the 2746-bus Polish network in 20 seconds on a standard laptop. 

\end{abstract}



\pagestyle{myheadings}
\thispagestyle{plain}
\markboth{Chance Constrained Optimal Power Flow}{D. Bienstock, M. Chertkov, S. Harnett}



The power grid can be considered one of the greatest engineering achievements of the 20th century, responsible for the economic well-being, social development, and resulting political stability of billions of people around the globe.  The grid is able to deliver on these goals with only occasional disruptions through significant control sophistication and careful long-term planning.

Nevertheless, the grid is under growing stress and the premise of secure electrical power delivered anywhere and at any time may become less certain.  Even though utilities have massively invested in infrastructure, grid failures, in the form of large-scale power outages, occur unpredictably and with increasing frequency.  In general, automatic grid control and regulatory statutes achieve robustness of operation as conditions display normal fluctuations, in particular approximately predicted inter-day trends in demand, or even unexpected single points of failure, such as the failure of a generator or tripping of a single line.  However, larger, unexpected disturbances can prove quite difficult to overcome.  This difficulty can be explained by the fact that automatic controls found in the grid are largely of an engineering nature (i.e. the flywheel-directed generator response used to handle short-term demand changes locally) and are largely \textit{not} of a data-driven, algorithmic and distributed nature. Instead, should an unusual condition arise, current grid operation relies on human input.  Additionally, only some real-time data is actually used by the grid to respond to evolving conditions.

All engineering fields can be expected to change as computing becomes ever more enmeshed into operations and massive amounts of real-time data become available.  In the case of the grid this change amounts to a challenge; namely how to migrate to a more algorithmic-driven grid in a cost-effective manner that is also seamless and gradual so as not to prove excessively disruptive -- because it would be impossible to rebuild the grid from scratch.   One of the benefits of the migration, in particular, concerns the effective integration of renewables into grids.  This issue is critical because large-scale introduction of renewables as a generation source brings with it the risk of large, random variability -- a condition that the current grid was not developed to accommodate.

\begin{figure}
\centering
\includegraphics[width=0.9\textwidth]{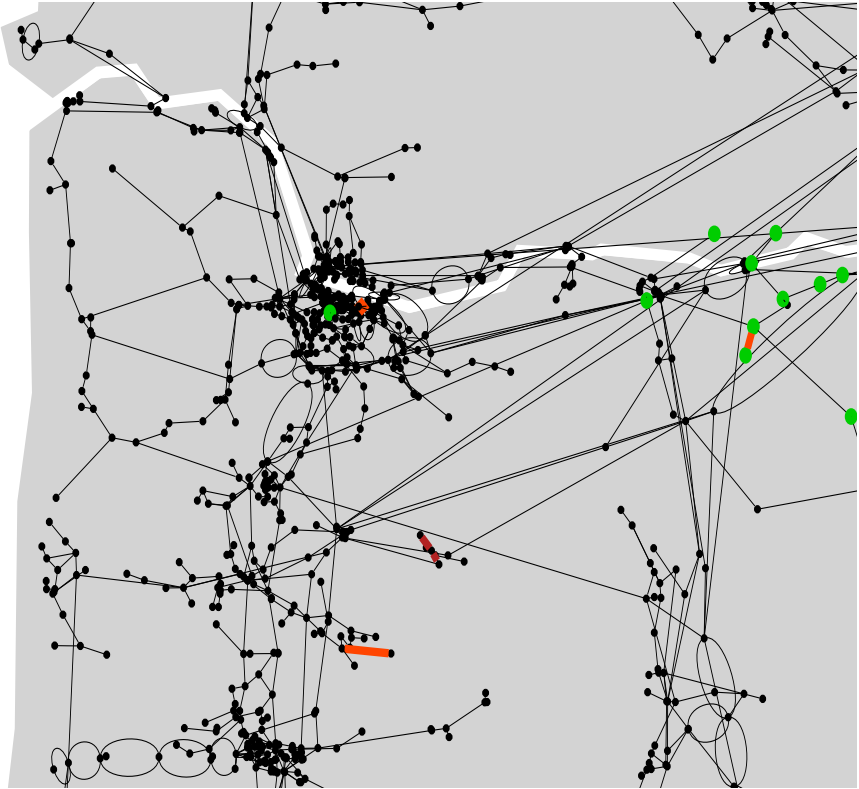}
\caption{
Bonneville Power Administration \cite{BPA}  shown in outline under 9\% wind penetration, where green dots mark actual wind farms. We set standard deviation to be 0.3 of the mean for each wind source. Our CC-OPF (with 1\% of overload set as allowable) resolved the case successfully (no overloads) and was computed in seconds, while the standard OPF showed 8 overloaded lines, all marked in color. Lines shown orange are at 4\% chance of overload. There are two dark red lines which are at 50\% of the overload while other (dark orange) lines show values of overload around 10\%.
 \label{BPA}} \end{figure}

This issue becomes clear when we consider how the grid sets generator output in ``real time''. This is typically performed at the start of every fifteen-minute (to an hour) period, or time window, using fixed estimates for conditions during the period.  More precisely, generators are dispatched so as to balance load (demand) and generator output at minimum cost, while adhering to operating limitations of the generators and transmission lines; estimates of the typical loads for the upcoming time window are employed in this computation.  This computational scheme, called Optimum Power Flow (OPF) or economic dispatch, can fail, dramatically, when renewables are part of the generation mix and (exogenous) fluctuations in renewable output become large.  By ``failure'' we mean, in particular, instances where a combination of generator and renewable outputs conspire to produce power flows that significantly exceed power line ratings.  When a line's rating is exceeded, the likelihood grows that the line will become tripped (be taken uncontrollably out of service) thus compromising integrity and stability of the grid.  If several key lines become tripped a grid would very likely become unstable and possibly experience a cascading failure, with large losses in serviced
demand.  This is not an idle assumption, since firm commitments to major renewable penetration are in place throughout the world. For example, 20\% renewable penetration by 2030 is a decree in the U.S. \cite{20-2030}, and similar plans are to be implemented in Europe, see e.g. discussions in \cite{CIGRE09,DENA,06Goz}. At the same time,  operational margins (between typical power flows and line ratings) are decreasing and expected to decrease.

A possible failure scenario is illustrated in Fig.~\ref{BPA} using as example
the U.S. Pacific Northwest regional grid data (2866 lines, 2209 buses, 176 generators and 18 wind sources),  where lines highlighted in red are jeopardized (flow becomes too high) with unacceptably high probability by fluctuating wind resources positioned along the Columbia river basin (green dots marking existing wind farms). We propose a solution that requires, as the only additional investment, accurate wind forecasts; but no change in machinery or significant operational procedures.  Instead, we propose an intelligent way to modify the optimization approach so as to mitigate risk; the approach is implementable as an efficient algorithm that solves realistic examples with thousands of buses and lines (such as the U.S. Pacific Northwest case) in a matter of seconds, and is thus only slightly slower than standard economic dispatch methods even on large-scale caes.

Maintaining line flows within their prescribed limits arises as a paramount operational criterion toward grid stability.  In the context of incorporating renewables into generation, a challenge emerges because a nominally safe way of operating a grid may become unsafe -- should an unpredictable (but persistent) change in renewable output occur, the resulting power flows may cause a line to persistently exceed its rating.  It is natural to assess the risk of such an event in terms of probabilities, because of the non-deterministic behavior of e.g. wind; thus in our proposed operational solution we will rely on techniques involving both mathematical optimization and risk analysis.

When considering a system under stochastic risk, an extremely large variety of events
that could pose danger might emerge.  Recent works \cite{10CPS,11CSPB,12BH} suggest that
focusing on instantons, or most-likely (dangerous) events, provides a
practicable route to risk control and assessment.  However, there may be
far too many comparably probable instantons, and furthermore, identifying such events does not answer the question of what to do about them.  In other words, we need a
computationally efficient methodology that not only identifies dangerous,
relatively probable events, but also mitigates them.

This paper suggests a new approach for handling the two challenges, that is to say, searching for the most probable realizations of line overloads under renewable generation, and correcting such situations through control actions, simultaneously and efficiently in one step. Our approach relies on methodologies recently developed in the optimization literature and known under the name of "Chance-Constrained" (CC) optimization \cite{06NS}.  Broadly speaking, CC optimization problems are optimization problems involving stochastic quantities, where constraints state that the probability of a certain random event is kept smaller than a target value.

To address these goals, we propose an enhancement of the standard OPF to be used in the economic dispatch of the controllable generators. We model each bus that houses a power source subject to randomness to include a random power injection, and reformulate the standard OPF in order to account for this uncertainty. The formulation minimizes the average cost of generation over the random power injections, while specifying a mechanism by which (standard, i.e. controllable) generators compensate in real-time for renewable power fluctuations; at the same time guaranteeing low probability that any line will exceed its rating. This last constraint is naturally formulated as a chance constraint -- we term out approach Chance-Constrained OPF, or CC-OPF.

This paper is organized as follows.  In Section \ref{sec:formulate} we motivate and present the various mathematical models used to describe how the grid operates, as well as our proposed methodology. We explain how to solve the models in Section \ref{sec:solve}. We then present, in Section \ref{sec:experiments} a number of examples to demonstrate the speed and usefulness of our approach. Section \ref{sec:discussions} summarizes the results and discusses the path forward.

\section{Formulating Chance-Constrained Optimum Power Flow Models}
\label{sec:formulate}

\subsection{Transmission Grids: Controls and Limits}

The power systems we consider in this paper are \textit{transmission grids} which operate at high voltages so as to convey power economically, with minimal losses, over large distances.  This is to be contrasted with \textit{distribution systems}; typically residential, lower voltage grids used to provide power to individual consumers.  From the point of view of wind-power generation, smooth operation of transmission systems is key since reliable wind sources are frequently located far away from consumption.

Transmission systems balance consumption/load and generation using a complex strategy that spans three different time scales (see e.g. \cite{00BV}).  At any point in time, generators produce power at a previously computed base level.  Power is generated (and transmitted) in the Alternating Current (AC) form. An essential ingredient toward stability of the overall grid is that all generators operate at a common frequency.  In real time, changes in loads are registered at generators through (opposite) changes in frequency.  A good example is that where there is an overall load increase.  In that case generators will marginally slow down -- frequency will start to drop. Then the so-called primary frequency control, normally implemented on gas and hydro plans with so-called ``governor" capability will react so as to stop frequency drift (large coal and nuclear units are normally kept on a time constant output). This is achieved by having each responding generator convey more power to the system, proportionally to the frequency change sensed. (In North America the proportionality coefficient is normally set to 5\% of the generator capacity for 0.5Hz deviation from the nominal frequency of 60Hz.) This reaction is swift and local, leading to stabilization of frequency across the system, however not necessarily at the nominal 60Hz value. The task of the secondary, or Automatic Gain Control (AGC), is to adjust the common frequency mismatch and thus to restore the overall balance between generation and consumption,  typically in a matter of minutes. Only some of the generators in a local area may be involved in this step. The final component in the strategy is the tertiary level of control, executed via the OPF algorithm,  typically run as frequently as every fifteen minutes (to one hour), and using estimates for loads during the next time window, where base (controllable) generator outputs are reset.  This is not an automatic step in the sense that a computation is performed to set these generator levels; the computation takes into account not only load levels but also other parameters of importance, such as line transmission levels.  Tertiary control computation, which is in the center of this paper, thus represents the shortest time scale where actual off-line and network wide (in contrast to automatic primary and secondary controls of frequency) optimal computations are employed. The three levels are not the only control actions used to operate a transmission system.  Advancing further in the time scale, OPF is followed by the so-called Unit Commitment (UC) computation, which schedules the switching on and off of large generation units on the scale of hours or even days.

A critical design consideration at each of the three control levels is that of maintaining ``stability'' of the grid.  The most important ingredient toward stable operation is \textit{synchrony} -- ultimately, all the generators of the network should stabilize thus locking, after a perturbation followed by a seconds-short transient, at the same frequency. Failure to do so not only proves inefficient but, worse, it threatens the integrity of the grid, ultimately forcing generators to shut down for protective reasons -- thus, potentially, causing a large, sudden change in power flow patterns (which may exceed equipment limits, see below) and possibly also an unrecoverable generation shortage.  A second stability goal is that of maintaining large voltages.  This is conducive to efficiency; lower voltage levels cause as a byproduct more generation (to meet the loads) and larger current values.  Not only is this combination inefficient, but in an extreme case it may make impossible to meet existing loads (so-called ``voltage collapse'' is a manifestation of this problem).  The third stability goal, from an operational perspective, is that of maintaining (line) power flows within established bounds. In long transmission lines, a large flow value will cause excessive voltage drop (an undesirable outcome as discussed). On a comparatively shorter line, an excessively large power flow across the line will increase the line temperature to the point that the line sags, and potentially arcs or trips due to a physical contact. For each line there is a given parameter, the \textit{line rating} (or \textit{limit}) which upper bounds flow level during satisfactory operation.

Of the three ``stability" criteria described above, the first two (maintaining synchrony and voltage) are a concern only in a truly nonlinear regime which under normal circumstances occur rarely. Thus, we focus on the third -- observing  line limits.

\subsection{OPF -- Standard Generation Dispatch (tertiary control)} \label{standardopf}

OPF is a key underlying algorithm in power engineering; see the review in \cite{Huneault1991} and e.g. \cite{94Kun,00BV}. The task of OPF, usually executed off-line at periodic intervals, is to reset generator output levels over a control area of the transmission grid, for example over the Bonneville Power Administration (BPA) grid shown in Fig.~\ref{BPA}. In order to describe OPF we will employ power engineering terms such ``bus'' to refer to a graph-theoretic vertex and ``line'' to refer to an edge. The set of all buses will be denoted by $\cV$, the set of lines, $\cE$ and the set of buses that house generators, $\cG$. We let $n = | \cV |$. A line joining buses $i$ and $j$ is denoted by $(i,j)$ indicating an arbitrary but fixed ordering. We assume that the underlying graph is connected, without loss of generality.\\

The generic OPF problem can be stated as follows:
\begin{itemize}
\item The goal is to determine the vector $p \in \R^{\cG}$, where for $i \in \cG$, $p_i$ is the output of generator $i$, so as to of minimize an objective function $c(p)$.  This function is, usually, a convex, separable quadratic function of $p$:
$$c(p) = \sum_{i \in \cG} c_i(p_i),$$
where each $c_i$ is convex quadratic.
\item The problem is endowed by three types of constraints: power flow, line
limit and generation bound constraints.
\end{itemize}

\vspace{.1in}
\noindent Among the constraints, the simplest are the generation bounds, which are box constraints on the individual $p_i$.  The thermal line limits place an upper bound on the power flows in each line.  We will return to these constraints later on; they are related to the physics model describing the transport of power in the given network, which is described by the (critical) power flow constraints. In the most general form, these are simply Kirchoff's circuit laws stated in terms of voltages (potentials) and power flows.  In this context, for each bus $i \in \cV$ its voltage $U_i$ is defined as $v_i e^{j\theta_i}$, where $v_i$ and $\theta_i$ are the voltage magnitude and phase angle at bus $i$.  Voltages can be used to derive expression for other physical quantities, such as current and, in particular, power, obtaining in the simplest case a system of quadratic equations on the voltage real and imaginary coordinates. See e.g. \cite{94Kun,00BV}.

The AC power flow equations  can constitute an obstacle to solvability of OPF (from a technical standpoint, they give rise to nonconvexities). In transmission system analysis a linearized version of the power flow equations is commonly used, the so-called ``DC-approximation''. In this approximation (a) all voltages are assumed fixed and re-scaled to unity; (b) phase differences between neighboring nodes are assumed small, $\forall (i,j)\in \cE:\quad |\theta_i-\theta_j|\ll 1$, (c) thermal losses are ignored (reactance dominates resistance for all lines). Then, the (real) power flow over line $(i,j)$, with line susceptance $\beta_{ij} \, (=\beta_{ji})$ is related linearly to the respective phase difference,
\begin{eqnarray}
f_{ij} & = & \beta_{ij}(\theta_i - \theta_j). \label{ohms}
\end{eqnarray}
Suppose, for convenience of notation, that we extend the vector $p$ to include an entry for every bus $i \in \cV$ with the proviso that $p_i = 0$ whenever $i \notin \cG$. Likewise, denote by $d \in \cR^{\cV}$ the vector of (possibly zero) demands  and by $\theta \in \cR^{\cV}$ the vector of phase angles.  Then, a vector $f$ of power flows is feasible if and only if
\begin{eqnarray}
\sum_{ij} f_{ij} & = & p_i - d_i, \ \ \ \mbox{for each bus $i$},
\end{eqnarray}
and in view of equation \eqref{ohms}, this can be restated as
\begin{eqnarray}
\theta_i \sum_{ij}\beta_{ij} \ \  - \ \ \sum_{ij}\beta_{ij}\theta_j  & = & p_i - d_i, \ \ \ \mbox{for each bus $i$}.
\end{eqnarray}
In matrix form this equation can be rewritten as follows:
\begin{eqnarray}
\quad && B \theta = p - d.   \label{pf1}
\end{eqnarray}
where the $n \times n$ matrix $B$ is a weighted-Laplacian defined as follows:
\begin{eqnarray}
    &&\forall i,j:\quad B_{ij}=\left\{\begin{array}{cc}
    -\beta_{ij}, & (i,j)\in\cE\\
    \sum_{k; (k,j)\in\cE}\beta_{kj},& i=j\\
    0,& \mbox{otherwise}\end{array}\right.,
    \label{B}\\
\end{eqnarray}
For future reference, we state some well known properties of Laplacians and the power flow system \eqref{pf1}.
\begin{LE} \label{trivial} The sum of rows of $B$ is zero and under the connectedness assumption for the underlying graph the rank of $B$ equals $n - 1$. Thus system \eqref{pf1} is feasible in $\theta$ if and only if
\begin{eqnarray}
\sum_i p_i & = & \sum_i d_i. \label{mybalance}
\end{eqnarray}
In other words: under the DC model the power flow Eqs.~\eqref{pf1} are feasible precisely when total generation equals total demand.   Moreover, if Eqs.~\eqref{pf1} are feasible, then for any index $1 \le j \le n$ there is a solution
with $\theta_j = 0$.
\QED \end{LE}\\

In summary, the standard DC-formulation OPF problem can be stated as the following constrained optimization problem:
\begin{eqnarray}
\mbox{{\bf OPF:}}\quad \min_p &&\, c(p),\quad \mbox{s.t.} \label{OPF-opt0}\\
\quad && B \theta = p - d \label{DC_flow1},\\
\quad && \forall i\in\cG:\quad p_i^{min}\leq p_i\leq p_i^{max},
\label{genconstr}
\label{GB}\\
\quad && \forall (i,j)\in\cE:\quad |f_{ij}|\leq f_{ij}^{max}, \label{tlconstr}
\label{TL}
\end{eqnarray}
Note that the $p_i^{min}, p_i^{max}$ quantities can be used to enforce the convention $p_i = 0$ for each $i \notin \cG$; if $i \in \cG$ then $p^{min}_i, p^{max}_i$ are lower and upper generation bounds which are generator-specific. Here, Constraint \eqref{TL} is the line limit constraint for $(i,j)$;  $f^{max}_{ij}$ represents the line limit (typically a thermal limit),  which is assumed to be strictly enforced in constraint \eqref{TL}. This conservative condition will be relaxed in the following.

\noindent Problem \eqref{OPF-opt0} is a convex quadratic program, easily solved using modern optimization tools. The vector $d$ of demands is \textit{fixed} in this problem and is obtained through estimation.  In practice, however, demand will fluctuate around $d$; generators then respond by adjusting their output (from the OPF-computed quantities) proportionally to the overall fluctuation as will be discussed below.

\noindent The scheme \eqref{OPF-opt0} works well in current practice, as demands do not substantially fluctuate on the time scale for which OPF applies. Thus the standard practice of solving \eqref{OPF-opt0} in the feasibility domain defined by Eq.~\eqref{DC_flow1},\eqref{B},\eqref{GB},\eqref{TL}, using demand forecasts based on historical data (and ignoring fluctuations) has produced a very reliable result - generation re-dispatch covering a span of fifteen minutes to an hour, depending on the system.

\subsection{Chance constrained OPF: motivation}

To motivate the problem we outline how generator output is modulated, in real time, in response to fluctuations of demand.  Suppose we have computed, using OPF, the output $p_i$ for each generator $i$ assuming constant demands $d$.  Let $\hat d(t)$ be the vector of real-time demands at time $t$.  Then so-called ``frequency control", or more properly, primary and secondary controls in combination that we will also call in the sequel ``affine" control, will reset generator outputs to quantities $\hat p_i(t)$ according to the following scheme
\begin{eqnarray}
\hat p_i(t) & = & p_i \ - \ \rho_i \sum_j (d_j - \hat d_j(t)) \ \ \ \mbox{for each $i \in \cG$}. \label{secondary}
\end{eqnarray}
In this equation, the quantities $\rho_i \ge 0$ are fixed and satisfy
$$ \sum_i \rho_i = 1.$$
Thus, from \eqref{secondary} we obtain
$$\sum_i \hat p_i(t) \ = \ \sum_i p_i - \sum_j (d_j - \hat d_j(t)) \ = \ \sum_j \hat d_j(t),$$
from Eq.~\eqref{mybalance}, in other words, demands are met. The quantities $\rho_i \ge 0$ are generator dependent but essentially chosen far in advance and without regard to short-term demand forecasts.

Thus, in effect, generator outputs are set in hierarchical fashion (first OPF, with adjustments as per \eqref{secondary} which is furthermore risk-unaware. This scheme has worked in the past because of the slow time scales of change in uncontrolled resources (mainly loads). That is to say, frequency control and load changes are well-separated. Note that a large error in the forecast or an under-estimation of possible $d$ for the next --e.g., fifteen minute-- period may lead to an operational problem in standard OPF (see e.g. the discussions in \cite{CAISO-2007,Makarov-wind}) because even though the vector $\hat p(t)$ is sufficient to meet demands, the phase angles $\hat \theta(t)$ computed from
$$ B \hat \theta(t) \ = \ \hat p(t) - \hat d(t)$$
give rise to real-time power flows
$$\hat f_{ij}(t) \ \doteq \ \beta_{ij}[\hat \theta_i(t) - \hat \theta_j(t)] $$
that violate constraints \eqref{tlconstr}. In fact, even the generator constraints \eqref{genconstr} may fail to hold. This has not been considered a handicap, however, simply because line trips due to overloading as a result of OPF-directed generator dispatch were (and still are) rare, primarily because the deviations $\hat d_i(t) - d_i$ will be small in the time scale of interest.  In effect, the risk-unaware approach that assumes constant demands in solving the OPF problem has worked well.

This perspective changes when renewable power sources such as wind are incorporated.  We assume that a subset $\cW$ of the buses holds uncertain power sources (wind farms); for each $j \in \cW$, write the amount of power generated by source $j$ at time $t$ as  $ \mu_j + \omega_j(t)$, where $\mu_j$ is the forecast output of farm $j$ in the time period of interest. For ease of exposition, we will assume in what follows that $\cG$ refers to the set of buses holding \textit{controllable} generators, i.e. $\cG \cap \cW = \emptyset$.  Renewable generation can be seamlessly incorporated into the OPF formulation \eqref{OPF-opt0}-\eqref{tlconstr} by simply setting $p_i = \mu_i$ for each $i \in \cW$. Assuming constant demands but fluctuating renewable generation, the application of the frequency control yields the following analogue to \eqref{secondary}:
\begin{eqnarray}
\hat p_i(t) & = & p_i \ - \ \rho_i \sum_{j \in \cW} \omega_j(t)\ \ \ \mbox{for each $i \in \cG$}, \label{windary}
\end{eqnarray}
e.g. if $\sum_{j \in \cW} \omega_j(t) > 0$, that is to say, there is a net increase in wind output, then (controllable) generator output will proportionally decrease.

Eq.~\eqref{windary} describes how generation will adjust to wind changes, under current power engineering practice.  The hazard embodied in this relationship is that the quantities $\omega_j(t)$ \textit{can} be large resulting in sudden and large changes in power flows, large enough to substantially overload power lines and thereby cause their tripping, a highly undesirable feature that compromises grid stability.  The risk of such overloads can be expected to increase (see \cite{CIGRE09}); this is due to a projected increase of renewable penetration in the future, accompanied by the decreasing gap between normal operation and limits set by line capacities. Lowering of the TL limits (the $f^{max}_{ij}$ quantities in Eqs.~\eqref{tlconstr}) can succeed in deterministically preventing overloads, but it also forces excessively conservative choices of the generation re-dispatch, potentially causing extreme volatility of the electricity markets. See e.g. the discussion in \cite{meyn-markets-volatility} on abnormal price fluctuations in markets that are heavily reliant on renewables.

\subsection{Using chance constraints}
\label{subsec:using-CC}

Power lines do not fail (i.e., trip) instantly when their flow thermal limits are exceeded. A line carrying flow that exceeds the line's thermal limit will gradually heat up and possibly sag, increasing the probability of an arc (short circuit) or even a contact with neighboring lines, with ground, with vegetation or some other object.  Each of these events will result in a trip. The precise process is extremely difficult to calibrate (this would require, among other factors, an accurate representation of wind strength and direction in the proximity of the line)\footnote{We refer the reader to \cite{usc} for discussions of line tripping during the 2003 Northeast U.S.-Canada cascading failure.}. Additionally, the rate at which a line overheats depends on its overload which may dynamically change (or even temporarily disappear) as flows adjust due to external factors; in our case fluctuations in renewable outputs. What \textit{can} be stated with certainty is that the longer a line stays \textit{overheated}, the higher the probability that it will trip -- to put it differently, if a line remains overheated long enough, then, after a span possibly measured in minutes, it will trip. In summary, (thermal) tripping of a line is primarily governed by the historical pattern of the overloads experienced by the line, and thus it may be influenced by the status of other lines (implicitly, via changes in power flows); further, exogenous factors can augment the impact of overloads.

Even though an exact representation of line tripping seems difficult, we can however state a practicable alternative. Ideally, we would use a constraint of the form ``{\it for each line, the fraction of the time that it exceeds its limit within a certain time window is small}''. Direct implementation of this constraint would require resolving dynamics of the grid over the generator dispatch time window of interest. To avoid this complication, we propose instead the following static proxy of this ideal model, a {\em chance constraint}: ``{\it we will require that the probability that a given line will exceed its limit is small}".

To formalize this notion, we assume: \\
\begin{itemize}
\item [{\bf W.1.}] For each $i \in \cW$, the (stochastic) amount of power generated by source $i$ is of the form $ \mu_i + \bm{\omega_i}$, where
\item [{\bf W.2.}] $\mu_i$ is constant, assumed known from the forecast, and $\bm{\omega_i}$ is a zero mean independent random variable with known standard deviation $\sigma_i$.
\end{itemize}
\vspace{.1in}

\noindent Here and in what follows, we use bold face to indicate uncertain quantities. Let $\bm{f_{ij}}$ be the flow on a given line $(i,j)$, and let $0 < \epsilon_{ij}$ be small. The chance constraint for line $(i,j)$, is:
\begin{eqnarray}
&& \quad P(\bm{f_{ij}} > f_{ij}^{max}) < \epsilon_{ij} \quad \mbox{and} \quad P(\bm{f_{ij}} < -f_{ij}^{max}) < \epsilon_{ij} \quad \quad \forall \ \ (i,j). \label{chance_constraint}
\end{eqnarray}
\noindent
One could alternatively use
\begin{eqnarray}
&& \quad P(|\bm{f_{ij}}| > f_{ij}^{max}) < \epsilon_{ij} \quad \forall \ \ (i,j),  \label{chance_constraint_joint}
\end{eqnarray}
which is less conservative than \eqref{chance_constraint}.  If \eqref{chance_constraint_joint} holds then so does \eqref{chance_constraint}, and if the latter holds then $P(|\bm{f_{ij}}| > f_{ij}^{max}) < 2 \epsilon_{ij}$. However, \eqref{chance_constraint} proves more tractable, and moreover we are interested in the regime where $\epsilon_{ij}$ is fairly small; thus we estimate that there is small practical difference between the two constraints; this will be verified by our numerical experiments.    Likewise, for a generator $g$ we will require that
\begin{eqnarray}
&& \quad P(\bm{p_g}  > p_{g}^{max}) < \epsilon_{g} \quad \mbox{and} \quad P(\bm{p_g} < p_{g}^{min}) < \epsilon_{g}.
\end{eqnarray}
The parameter $\epsilon_g$ will be chosen \textit{extremely} small, so that for all practical purposes the generator's will be guaranteed to stay within its bounds.

Chance constraints \cite{chance}, \cite{charnes}, \cite{millerwagner} are but one possible methodology for handling uncertain data in optimization. Broadly speaking, this methodology fits within the general field of stochastic optimization.  Constraint \eqref{chance_constraint} can be viewed as a ``value-at-risk'' statement; the closely-related ``conditional value at risk'' concept provides a (convex) alternative, which roughly stated constrains the expected overload of a line to remain small, \textit{conditional} on there being an overload  (see \cite{06NS} for definitions and details).   Even though alternate models are
possible, we would would argue that our chance constrained approach is reasonable (in fact: compelling) in view of the nature of the line tripping process we discussed above.

\cite{Zhang2011} considers the standard OPF problem under stochastic demands, and describes a method that computes \textit{fixed} generator output levels to be used throughout the period of interest, independent of demand levels.  In order to handle variations in demand, \cite{Zhang2011} instead relies on the concept of a \textit{slack bus}. A slack bus is a fixed node that is assumed to compensate for all generation/demand mismatches -- when demand exceeds generation the slack bus injects the shortfall, and when demand is smaller than generation the slack bus absorbs the generation excess.  A vector of generations is acceptable if the probability that each system component operates within acceptable bounds is high -- this is a chance constraint.  To tackle this problem \cite{Zhang2011} proposes a simulation-based local optimization system consisting of an outer loop used to assess the validity of a control (and estimate its gradient) together with an inner loop that seeks to improve the control. Experiments are presented using a 5-bus and a 30-bus example.  The approach in \cite{Zhang2011}, though universal (i.e., applicable to any type of exogenous distribution),
requires a number of technical assumptions and elaborations to guarantee convergence and feasibility and appears to entail a very high computational cost.

Another related study \cite{12VMLA} describes a security constrained optimization for reserve scheduling with fluctuating wind generation. Security constraints are represented via a set of outage scenarios. Finite time optimization horizon (24 steps, one per hour) is considered. As in our study, the standard DC power flow approach is adopted and no load uncertainty is considered. 
However, in a contrast with our setting wind generation in \cite{12VMLA} is assumed localized (located at a single bus) in a small, 30 bus network. Under the 
assumption that the probability distribution function of the aggregated wind resources is known, the cost of generation over the time horizon is optimized under the chance constrains for line overloads and generation limits. The stochastic optimization is tackled in \cite{12VMLA} via transformation to a convex optimization on the expense of a number of approximations, including scenario approach to the wind forecast, look-up table modeling of the security constraints, and a heuristic scheme for a bilinear descent iteration, with no convergence to global optimum of the nonlinear problem guaranteed.

Chance constrained optimization has also been discussed recently in relation to the Unit Commitment problem, which concerns discrete-time planning for operation of large generation units on the scale of hours-to-months, so as to account for the long-term wind-farm generation uncertainty \cite{Ozturk2004,08WSL,11ZSHY}.

\subsection{Uncertain power sources}
\label{uncertainwind}

The physical assumptions behind our model of uncertainty are as follows. Independence of fluctuations at different sites is due to the fact that the wind farms are sufficiently far away from each other. For the typical OPF time span of 15 min and typical wind speed of $10m/s$, fluctuations of wind at the farms more than $10km$ apart are not correlated. We also rely on the assumption that transformations from wind to power at different wind farms is not correlated.

To formulate and calibrate our models, we make simplifying assumptions that are approximately consistent with our general physics understanding of fluctuations in atmospheric turbulence; in particular we assume Gaussianity of $\bm{\omega_i}$ \footnote{Correlations of velocity within the correlation  time of 15 min, roughly equivalent to the time span between the two consecutive OPF, are approximately Gaussian. The assumption is not perfect, in particular because it ignores significant up and down ramps possibly extending tails of the distribution in the regime of really large deviations.}. We will also assume that only a standard weather forecast (coarse-grained on minutes and kilometers) is available, and no systematic spillage of wind in its transformation to power is applied\footnote{See \cite{CIGRE09}, for some empirical validation.}.

There is an additional and purely computational reason for the Gaussian assumption, which is that under this assumption chance constraint \eqref{chance_constraint} can be captured using
in an optimization framework that proves particularly computationally practicable.
We will also consider a data-robust version of our chance-constrained problem where the
parameters for the Gaussian distributions are assumed unknown, but lying in a window.  This
allows both for parameter mis-estimation and for model error, that is to say the
implicit approximation of non-Gaussian distributions with Gaussians; our approach, detailed in Section \ref{subsec:ambiguous}, remains computationally sound in this robust setting.

Other fitting distributions considered in the wind-modeling literature,  e.g. Weibull distributions and logistic distributions \cite{08BDL,nrelweibull}, will be discussed later in the text as well \footnote{Note that the fitting approach of \cite{08BDL,nrelweibull} does not differentiate between typical and atypical events and assumes that the main body and the tail should be modeled using a simple distribution  with only one or two fitting parameters.  Generally this assumption is not justified as the physical origin of the typical and anomalous contributions of the wind, contributing to the main body and the tail of the distribution respectively, are rather different. Gaussian fit (of the tail)  -- or more accurately, faster than exponential decay of probability in the tail for relatively short-time (under one hour) forecast -- would be reasonably consistent with phenomenological modeling of turbulence generating these fluctuations.}. In particular, we will demonstrate on out-of-sample tests that the computationally advantageous Gaussian modeling of uncertainty allows as well to model effects of other distributions.

\subsection{Affine Control}
\label{subsec:affine}

Since the power injections at each bus are fluctuating, we need a control to ensure that generation is equal to demand at all times within the time interval between two consecutive OPFs. We term the joint result of the primary frequency control and secondary frequency control the \textit{affine} control. The term will intrinsically assume that all governors involved in the controls respond to fluctuations in the generalized load (actual demand which is assumed frozen minus stochastic wind resources) in a proportional way, however with possibly different proportionality coefficients. Then, the random variable ${\bf \omega}$ dependent version of Eq.~(\ref{windary}) becomes
\begin{eqnarray}
\forall \, \mbox{bus} \, i\in\cG:\quad \bm{p_i} = \bar p_i - \alpha_i \sum_{j \in \cW} \bm{\omega_j}.\quad
\alpha_i \ge 0.
\label{affinecontrol}
\end{eqnarray}
Here the quantities $\bar p_i \ge 0$ and $\alpha_i \ge 0$ are design variables satisfying (among other constraints) $\sum_{i \in \cG} \alpha_i = 1$. Notice that we do not set any $\alpha_i$ to a standard (fixed) value, but instead leave the optimization to decide the optimal value. (In some cases it may even be advantageous to allow negative $\alpha_i$ but we decided not to consider such a drastic change of current policy in this study.) The generator output $p_i$ combines a fixed term $\bar p_i$ and a term which varies with wind, $-\alpha_i \sum_{j \in \cW} \bm{\omega_j}$. Observe that $\sum_i p_i = \sum_i \bar p_i -\sum_i \bm{\omega_j}$, that is, the total power generated equals the average production of the generators minus any additional wind power above the average case.

This affine control scheme creates the possibility of requiring a generator to produce power beyond its limits. With unbounded wind, this possibility is inevitable, though we can restrict it to occur with arbitrarily small probability, which we will do with  additional chance constraints for all controllable generators, $\forall g\in\cG$,
\begin{align}\label{chance_constraint_gen}
    \quad P(p_g^{min} \leq \bar p_g - \alpha_g
    \sum_{j \in \cW} \bm w_j \leq p_g^{max}) > 1 - \epsilon_g.
\end{align}

\subsection{CC-OPF: Brief Discussion of Solution Methodology}
\label{subsec:brief}

Our methodology applies and develops general ideas of \cite{06NS} to the power engineering setting of generation re-dispatch under uncertainty. In Section \ref{subsec:formal} we will
provide a generic formulation of our chance-constrained OPF problem that is valid under
the assumption of linear power flow laws and statistical independence of wind fluctuations
at different sites, while using control law \eqref{affinecontrol} to specify standard
generation response to wind fluctuations.  Under the additional assumption of
Gaussianity, this generic formulation is reduced to a specific deterministic optimization problem of Section \ref{subsec:conic}.  Moreover this deterministic optimization over $\overline{p}$ and $\alpha$ is a convex optimization problem, more precisely, a Second-Order Cone Program (SOCP) \cite{BVbook,SOCP}, allowing an efficient computational implementation discussed in detail in Section \ref{subsec:cuttingplane}. We will term this SOCP, which assumes point
estimates for the wind distributions, the \textit{nominal} problem. As indicated
above, we also discuss how the SOCP formulation can be extended to account for data-related uncertainty in the parameters of the Gaussian distributions in Section \ref{subsec:ambiguous},
obtaining a \textit{robust} version of the chance-constrained optimization problem.

Let us emphasize that many of our assumptions leading to the computational efficient nominal formulation are not restrictive and allow natural generalizations. In particular, using techniques from \cite{06NS}, it is possible to relax the phenomenologically reasonable but approximately validated assumption of wind source Gaussianity (validated according to actual measurements of wind, see \cite{08BDL,nrelweibull} and references therein). For example, using only the mean and variance of output at each wind farm, one can use Chebyshev's inequality to obtain a similar though more conservative formulation. And following \cite{06NS} we can also obtain convex approximations to~\eqref{chance_constraint} which are tighter than Chebyshev's inequality, for a large number of empirical distributions discussed in the literature. The
data-robust version of our algorithm provides a methodologically sound (and
computationally efficient) means to protect against data and model errors; moreover we will perform (below) out-of-sample experiments involving the controls computed with the nominal approach; first to investigate the effect of parameter estimation errors in the Gaussian case, and, second, to gauge the impact of non-Gaussian wind distributions.

\section{Solving the Models}
\label{sec:solve}

\subsection{Chance-constrained optimal power flow: formal expression}
\label{subsec:formal}

Following the W.1 and W.2 notations, Eqs.~\eqref{affinecontrol} explain the \textit{affine} control, given that the $\alpha_i$ are decision variables in our CC-OPF, additional to the standard $\bar p_i$ decision variables already used in the standard OPF (\ref{OPF-opt0}). For $i \notin \cW$ write $\mu_i = 0$, thereby obtaining a vector $\mu \in \R^n$. Likewise, extend $\bar p$ and $\alpha$ to vectors in $\R^n$ by writing $\bar p_i = \alpha_i = 0$ whenever $i \notin \cG$.  \\

\noindent {\bf Definition.}  We say that the pair $\bar p, \alpha$ is \textit{viable} if the generator outputs under control law \eqref{affinecontrol}, together with the uncertain outputs, always exactly match total demand.  \\

\noindent The following simple result characterizes this condition as well as other basic properties of the affine control. Here and below, $e \in \R^n$ is the vector of all 1's. \\

\noindent \begin{LE}\label{easy1} Under the control law \eqref{affinecontrol} the net output of bus $i$ equals
\begin{eqnarray}
&& \bar p_i + \mu_i  - d_i + \bm{\omega} - \alpha_i (e^T \bm{\omega}), \label{outi}
\end{eqnarray}
and thus the (stochastic) power flow equations can be written as
\begin{eqnarray}
 B \bm{\theta} & = & \bar p + \mu - d + \bm{\omega} - (e^T \bm{\omega}) \alpha.
\label{stochthetasol}
\end{eqnarray}
Consequently, the pair $\bar p, \alpha$ is viable if and only if
\begin{eqnarray}
&& \sum_{i\in\cV} ( \bar p _i +      \mu_i - d_i) = 0. \label{bal2}
\end{eqnarray}
\end{LE}
\noindent {\em Proof.}  Eq.~\eqref{outi} follows by definition of the $\bar p$, $\mu$, $d$ vectors and the control law.  Thus Eq.~\eqref{stochthetasol} holds. By Lemma \ref{trivial} from Eq.~\eqref{stochthetasol} one gets that $\bar p, \alpha$ is viable iff
\begin{eqnarray}
0 & = & \sum_{i = 1}^n (\bar p_i - (e^T \bm{\omega}) \alpha_i + \mu_i + \bm{\omega_i} - d_i) \nonumber \\
&=& \sum_{i} (\bar p_i + \mu_i - d_i), \label{rbalance3}
\end{eqnarray}
since by construction $\sum_i \alpha_i = 1$.
\QED \\

\noindent {\bf Remark.} Equation \eqref{bal2} can be interpreted as stating the condition that expected total generation must equal total demand; however the Lemma contains a rigorous proof of this fact.\\

As remarked before, any $(n-1) \times (n-1)$ matrix obtained by striking out the same column and row of $B$ is invertible.  For convenience of notation we will assume that bus $n$ is neither a generator nor a wind farm bus, that is to say, $n \notin \cG \cup \cW$, and we denote by $\hat B$ the submatrix obtained from $B$ by removing row and column $n$, and write
\begin{eqnarray}
\breve \cB & = & \left( \begin{array}{c c}
\hat B^{-1} & 0\\
0 & 0 \end{array} \right). \label{brevedef}
\end{eqnarray}
Further, by Lemma \ref{trivial} we can assume without loss of generality that $\bm{\theta_n} = 0$. The following simple result will be used in the sequel.

\noindent \begin{LE}\label{easy} Suppose the pair $\bar p, \alpha$ is viable. Then under the control law \eqref{affinecontrol} a vector of (stochastic) phase angles is
\begin{eqnarray}
\bm{\theta} & = & \bar \theta + \breve \cB (\bm{\omega} - (e^T \bm{\omega}) \alpha ), \label{stochtheta} \ \ \mbox{where} \\
\bar \theta  & = & \breve \cB (\bar p  + \mu - d). \label{barthetadef}
\end{eqnarray}
As a consequence,
\begin{eqnarray}
 \mathbb{E_{\bm{\omega}}} \bm{\theta} & = & \bar \theta, \label{exptheta}
\end{eqnarray}
and given any line $(i,j)$,
\begin{eqnarray}
\mathbb{E_{\bm{\omega}}} \bm{f_{ij}} & = & \beta_{ij}(\bar \theta_i - \bar \theta_j). \label{expflow}
\end{eqnarray}
Furthermore, each quantity \bmm{\theta_i} or \bmm{f_{ij}} is an affine function of the random variables \bmm{\omega_i}.
\end{LE}\\
\noindent {\em Proof.}  For convenience we rewrite system \eqref{stochthetasol}: $B \bm{\theta} \, = \, \bar p + \mu - d + \bm{\omega} - (e^T \bm{\omega}) \alpha.$ Since $\bar p, \alpha$ is viable, this system is always feasible, and since the sum of rows of $B$ is zero, its last row is redundant. Therefore Eq.~\eqref{stochtheta} follows since $\bm{\theta_n} = 0$, and Eq.~\eqref{exptheta} holds since $\bm{\omega}$ has zero mean. Since $\bm{f_{ij}}  =  \beta_{ij}  (\bm{\theta_i} - \bm{\theta_j})$ for all $(i,j),$ Eq.~\eqref{expflow} holds. From this fact and  Eq.~\eqref{stochtheta} it follows that $\bm{\theta}$ and $\bm{f}$ are affine functions of $\bm{\omega}$.
\QED \\

\noindent Using this result we can now give an initial formulation to our chance-constrained problem; with discussion following.
\begin{eqnarray}
&& \mbox{{\bf CC-OPF:}}\quad\min \ \mathbb{E}_{\bm{\omega}}\left[c(\bar p - (e^T \bm{\omega})\alpha ) \,  \right]  \label{CC-OPF}\\
&&        \mbox{s.t. } \ \
 \sum\limits_{i\in\cG} \alpha_i  =  1,\quad  \alpha \ge 0,  \quad \bar p \ge 0\label{conic-first}\\
&& \sum\limits_{i\in\cV} ( \bar p _i +        \mu_i - d_i) = 0 \label{thetabar2} \\
&& \quad \quad B \bar \theta  \ = \ \bar p  + \mu - d, \label{thetabar}\\
&& \nonumber \\
&& \mbox{for all lines $(i,j)$:} \nonumber\\
&& \quad P \left( \beta_{ij} (\bar \theta_i - \bar \theta_j + [ \breve \cB (\bm{\omega} - (e^T \bm{\omega}) \alpha) ]_i - [\breve \cB (\bm{\omega} - (e^T \bm{\omega}) \alpha)]_j \, ) > f_{ij}^{max}\right)  <  \epsilon_{ij} \nonumber \\
&& \label{longchance1} \\
&& \quad P \left( \beta_{ij} (\bar \theta_i - \bar \theta_j + [ \breve \cB (\bm{\omega} - (e^T \bm{\omega}) \alpha) ]_i - [\breve \cB (\bm{\omega} - (e^T \bm{\omega}) \alpha)]_j \, ) < -f_{ij}^{max}\right)  <  \epsilon_{ij} \nonumber \\
&& \label{longchance2}\\
&& \mbox{for all generators $g$:} \nonumber\\
&& \quad P \left(\bar p_g - (e^T\bm{\omega}) \alpha_i > p^{max}_g \right) < \epsilon_g \ \ \mbox{and} \ \ P \left(\bar p_g - (e^T\bm{\omega}) \alpha_i < p^{min}_g \right)  < \epsilon_g. \label{genchance}
\end{eqnarray}
The variables in this formulation are $\bar p, \alpha$ and $\bar \theta$. Constraint \eqref{conic-first} simply states basic conditions needed by the affine control.   Constraint \eqref{thetabar2} is \eqref{bal2}.  Constraints \eqref{thetabar}, \eqref{longchance1} and \eqref{longchance2} express our chance constraint, in view of Lemma \ref{easy}.

The objective function is the expected cost incurred by the stochastic generation vector
$$ \bm{p} =  \bar p - (e^T \bm{\omega})\alpha$$
over the varying wind power output $\bm w$. In standard power engineering practice generation cost is convex, quadratic and separable, i.e. for any vector $p$, $c(p) = \sum_i c_i(p_i)$ where each $c_i$ is convex quadratic. Note that for any $i \in \cG$ we have
$$ \bm{p_i}^2 \ = \ \bar p_i^2 + (e^T \bm{\omega})^2 \alpha_i^2 - 2 e^T \bm{\omega} \bar p_i \alpha_i,$$
from which we obtain, since the $\omega_i$ have zero mean,
$$ \mathbb{E}_{\bm w} (\bm{p_i}^2) \ = \ \bar p_i^2 + var(\bm{\Omega}) \alpha_i^2,$$
where ``var'' denotes variance and $\bm{\Omega} \doteq \sum_j \bm{\omega_j}$. It follows that the objective function can be written as
\begin{eqnarray}
\mathbb{E_{\bm{\omega}}} c(\bm{p}) & = & \sum_{i \in \cG} \left\{ c_{i1}\left(\bar p_i^2 + var(\bm{\Omega}) \alpha_i^2 \right) + c_{i2} \bar p_i + c_{i3} \right\}
\end{eqnarray}
where $c_{i1} \ge 0$ for all $i \in \cG$.  Consequently the objective function is convex quadratic, as a function of $\bar p$ and $\alpha$.

The above formulation is the formal statement for our optimization problem. Even though its objective is convex in cases of interest, the formulation is not in a form that can be readily exploited by standard optimization algorithms. Below we will provide an efficient approach to solve relevant classes of problems with the above form; prior to that we need a technical result.  We will employ the following notation:
\begin{itemize}
\item For $j \in \cW$, the variance of $\bm{\omega_j}$ is denoted by $\sigma_j^2$.
\item For $1 \le i, j \le n$ let $\pi_{ij}$ denote the $i, j$ entry of
the matrix $\breve \cB$ given above, that is to say,
\begin{eqnarray}
\pi_{ij} \ = \ \left\{\!\!
\begin{array} {cc} [\hat{B}^{-1}]_{ij},\! &  i < n,  \\
0,\!\!\! & \mbox{otherwise.}
		\end{array}
	\right.\!\!\!,\
\end{eqnarray}
\item Given $\alpha$, for $1 \le i \le n$ write
\begin{eqnarray}
\delta_i & \doteq & [\breve \cB \alpha]_i \ = \  \left\{\!\! \begin{array} {ll}
		[\hat B^{-1} \alpha]_i,\!\!\!  & i < n,\\
		0,\!\!\! & \mbox{otherwise.}  \label{dnumber}
		\end{array}
	\right.
\end{eqnarray}
\end{itemize}

\begin{LE}\label{independent}  Assume
that the $\bm{\omega_i}$ are
independent random variables. Given $\alpha$, for any line $(i,j)$,
\begin{eqnarray}
var( \bmm{f_{ij}} )& = & \beta_{ij}^2 \sum_{k \in \cW} \sigma_k^2 (\pi_{i k} - \pi_{jk} - \delta_i + \delta_j)^2.  \label{ole'}
\end{eqnarray}
\end{LE}
\noindent {\em Proof.}  Using $\bm{f_{ij}}  =  \beta_{ij} [ \bm{\theta_i} - \bm{\theta_j}]$ and eq. \eqref{stochtheta} we have that
\begin{eqnarray}
\bm{f_{ij}} - \mathbb{E_{\bm{\omega}}} \bm{f_{ij}} & = & \beta_{ij} ( [\breve \cB(\bm{\omega} - (e^T\bm{\omega})\alpha) ]_i - [\breve \cB(\bm{\omega} - (e^T\bm{\omega})\alpha) ]_j) \ = \ \\
& = & \beta_{ij} \left( [\breve \cB \bm{\omega}]_i - [\breve \cB \bm{\omega}]_j - (e^T \bm{\omega}) \delta_i + (e^T \bm{\omega}) \delta_j \right) \ = \ \\
& = & \beta_{ij} \sum_{k \in \cW} (\pi_{ik} - \pi_{jk} - \delta_i + \delta_j) \bm{\omega_k},
\end{eqnarray}
since by convention $\bm{\omega_i} = 0$ for any $i \notin \cW$.  The result now follows. \QED\\

\noindent {\bf Remark.} Lemma \ref{independent} holds for any distribution of the\bmm{\omega_i} so long as independence is assumed.  Similar results are easily obtained for higher-order moments of the\bmm{f_{ij}}.

\subsection{Formulating the chance-constrained problem as a conic program}
\label{subsec:conic}

In deriving the above formulation \eqref{CC-OPF}-\eqref{genchance} for CC-OPF we assumed that the $\bm{\omega_i}$ random variables have zero mean.  To obtain an efficient solution procedure we will additionally assume that they are (a) pairwise independent and (b) normally distributed. These assumptions were justified in Section \ref{uncertainwind}. Under the assumptions, however, since the\bmm{f_{ij}} are affine functions of the $\bm{\omega_i}$ (because the $\bm{\theta}$ are, by eq. \eqref{stochtheta}), it turns out that there is a simple restatement of the chance-constraints \eqref{longchance1},
\eqref{longchance2} and \eqref{genchance} in a computationally practicable form.  See \cite{06NS} for a general treatment of linear inequalities with stochastic
coefficients.  For any real $0 < r < 1$ we write $\eta(r) = \phi^{-1}(1 - r)$, where
$\phi$ is the cdf of a standard normally distributed random variable.

\begin{LE}\label{normalline}
Let $\bar p$, $\alpha$ be viable.  Assume that the $\bm{\omega_i}$ are normally distributed and independent. Then: \\

\noindent For any line $(i,j)$, $P(\bm{f_{ij}} > f^{max}_{ij}) < \epsilon_{ij}$ and $P(-\bm{f_{ij}} > f^{max}_{ij}) < \epsilon_{ij}$ if and only if
\begin{eqnarray}
&&  \beta_{ij}|\bar \theta_i - \bar \theta_j| \ \leq \ f^{max}_{ij}-\eta(\epsilon_{ij}) \left[ \beta_{ij}^2 \sum_{k \in \cW} \sigma_k^2 (\pi_{i k} - \pi_{jk} - \delta_i + \delta_j)^2\right]^{1/2}  \label{up-down}
\end{eqnarray}
where as before $\bar \theta =  \breve \cB (\bar p  + \mu - d)$ and
$\delta = \breve \cB \alpha$.  \\

\noindent For any generator $g$,
$P \left(\bar p_g - (e^T\bm{\omega}) \alpha_i > p^{max}_g \right) < \epsilon_g$ and $P \left(\bar p_g - (e^T\bm{\omega}) \alpha_i < p^{min}_g \right) < \epsilon_g$ iff
\begin{eqnarray}
p^{min}_g + \eta(\epsilon_g) \left( \sum_{k \in \cW}\sigma^2_k\right)^{1/2} \ \le \ \bar p_g \ \le \ p^{max}_g - \eta(\epsilon_g) \left( \sum_{k \in \cW}\sigma^2_k\right)^{1/2}. \label{genchanceud}
\end{eqnarray}
\end{LE}\\
\noindent {\em Proof.} By Lemma \ref{easy}, $\bm{f_{ij}}$ is an affine function
of the $\bm{\omega_i}$; under the assumption it follows that $\bm{f_{ij}}$
is itself normally distributed. Thus, $P(\bm{f_{ij}} > f^{max}_{ij}) < \epsilon_{ij}$ iff
\begin{eqnarray}
 \mathbb{E_{\bm{\omega}}} \bm{f_{ij}} + \eta(\epsilon_{ij}) \, var(\bm{f_{ij}})  & \le & f^{max}_{ij}, \label{upcase}
\end{eqnarray}
and similarly, $P(\bm{f_{ij}} < -f^{max}_{ij}) < \epsilon_{ij}$ iff
\begin{eqnarray}
 \mathbb{E_{\bm{\omega}}} \bm{f_{ij}} - \eta(\epsilon_{ij}) \, var(\bm{f_{ij}})  & \ge & -f^{max}_{ij}. \label{downcase}
\end{eqnarray}
Lemma \ref{easy} gives $ \mathbb{E_{\bm{\omega}}} \bm{f_{ij}}  =  \beta_{ij}(\bar \theta_i - \bar \theta_j)$ while by Lemma \ref{independent},
$var( \bmm{f_{ij}} ) =  \beta_{ij}^2 \sum_{k \in \cW} \sigma_k^2 (\pi_{i k} - \pi_{jk} - \delta_i + \delta_j)^2$.  Substituting these values into \eqref{upcase} and
\eqref{downcase} yields \eqref{up-down}.   The proof of \eqref{genchanceud} is similar.\QED \\

\noindent {\bf Remarks} Eq. \eqref{up-down} highlights the difference between e.g.
our chance constraint for lines, which requires that $P(\bm{f_{ij}} > f^{max}_{ij}) < \epsilon_{ij}$ and that $P(\bm{f_{ij}} < -f^{max}_{ij}) < \epsilon_{ij}$, and the stricter
requirement that $P(|\bm{f_{ij}}| > f^{max}_{ij}) < \epsilon_{ij}$ which amounts to
\begin{eqnarray}
&& P(\bm{f_{ij}} > f^{max}_{ij}) + P(\bm{f_{ij}} < -f^{max}_{ij}) < \epsilon_{ij}.\label{strictupdown}
\end{eqnarray}
Unlike our requirement, which is captured by \eqref{up-down},
the stricter condition \eqref{strictupdown} does not admit a compact statement.\\

We can now present a formulation of our chance-constrained optimization as a
convex optimization problem, on variables $\bar p, \alpha, \bar \theta, \delta$
and $s$:
\begin{eqnarray}
\min \, \sum_{i \in \cG} \left\{ c_{i1} \bar p_i^2 + \left(\sum_k \sigma^2_k\right) c_{i1} \alpha_i^2 + c_{i2} \bar p_i + c_{i3} \right\} \quad && \ \label{obj}\\
\mbox{for $1 \le i \le n-1$:}\quad \quad \quad \ \quad \quad \sum_{j = 1}^{n-1} \hat B_{ij} \, \delta_j = \alpha_i\quad
&& \ \ \label{deltadef} \\ 
\mbox{for $1 \le i \le n-1$:} \quad \sum_{j = 1}^{n-1} \hat B_{ij} \bar \theta_j \ - \bar p_i   =  \mu_i - d_i \quad
 &&  \label{Edef}\\
&& \nonumber \\
\sum_i \alpha_i = 1, \quad \alpha \ge 0, \quad \bar p \ge 0 \quad && \label{basic}\\
\bar p_n = \alpha_n = \delta_n  =  \bar \theta_n = 0 \quad && \\
&& \nonumber \\
  \beta_{ij} |\bar \theta_i - \bar \theta_j| + \beta_{ij} \eta(\epsilon_{ij}) \, s_{ij} \ \leq \ f^{max}_{ij}  \quad && \ \forall \ (i,j)\quad \label{sdef}\\
&& \nonumber \\
  \left[\sum_{k \in \cW} \sigma_k^2 (\pi_{i k} - \pi_{jk} - \delta_i + \delta_j)^2\right]^{1/2} \ - \ s_{ij} \ \leq \ 0  \quad && \ \forall \ (i,j)\quad\label{varf}\\
&& \nonumber \\
-\bar p_g + \eta(\epsilon_g) \left( \sum_{k \in \cW}\sigma^2_k\right)^{1/2} \, \le \, -p^{min}_g   \quad && \ \forall \ g \in \cG \label{gen1}\\
\bar p_g  + \eta(\epsilon_g) \left( \sum_{k \in \cW}\sigma^2_k\right)^{1/2}\, \le \ p^{max}_g \quad && \ \forall \ g \in \cG \label{gen2}
\end{eqnarray}
In this formulation, the variables $s_{ij}$ are auxiliary and introduced to facilitate the discussion below -- since $\eta_{ij} \ge 0$ without loss of generality \eqref{varf} will hold as an equality. Constraints \eqref{varf}, \eqref{gen1} and \eqref{gen2} are \textit{second-order cone} inequalities \cite{BVbook}.  A problem of the above form is solvable in polynomial time using well-known methods of convex optimization; several commercial software tools such as Cplex \cite{CPLEX}, Gurobi \cite{gurobi}, Mosek \cite{mosek} and others are available. Constraint \eqref{deltadef} is equivalent to $\delta_i = [\hat B^{-1} \alpha]_i$ (as we did in \eqref{dnumber}), however $\hat B^{-1}$ can be seen to be totally dense, whereas $\hat B$ is very sparse for typical grids.  Constraint \eqref{varf} can be relatively dense -- the sum has a term for each farm.  However as a percentage of the total number of buses this can be expected to be small.

\subsection{Solving the conic program}
\label{subsec:cuttingplane}

Even though optimization theory guarantees that the above problem is efficiently solvable, experimental testing shows that in the case of large grids (thousands of lines) the problem proves challenging.  For example, in the Polish 2003-2004 winter peak case\footnote{Available with MATPOWER \cite{MATPOWER}}, we have $2746$ buses, $3514$ lines and $8$ wind farms, and Cplex \cite{CPLEX} reports (after pre-solving) $36625$ variables and $38507$ constraints, of which $6242$ are conic. On this problem, a recent version of Cplex on a (current) 8-core workstation ran for $3392$ seconds (on 16 parallel threads, making use of ``hyperthreading'') and was unable to produce a feasible solution.  On the same problem Gurobi reported ``numerical trouble'' after 31.1 cpu seconds, and stopped.

In fact, all of the commercial solvers \cite{CPLEX,gurobi,mosek} we experimented with reported numerical difficulties with problems of this size.  Anecdotal evidence indicates that the primary cause for these difficulties is not simply the size, but also to a large degree \textit{numerics} in particular poor conditioning due to the entries in the matrix $B$. These are susceptances, which are inverses of reactances, and often take values in a wide range.

To address this issue we implemented an effective algorithm for solving problem \eqref{obj}-\eqref{gen2}.  For brevity we will focus on constraints \eqref{varf} (\eqref{gen1} and \eqref{gen2} are similarly handled).  For a line $(i,j)$ define
$$C_{ij}(\delta) \doteq \left( \sum_{k \in \cW} \sigma_k^2(\pi_{i k} - \pi_{jk} - \delta_i + \delta_j)^2 \right)^{1/2}.$$
Constraint \eqref{varf} can thus be written as $ C_{ij}(\delta) \le s_{ij}$. For completeness, we state the following result:
\begin{LE} \label{outer} Constraint \eqref{varf} is equivalent to the infinite
set of linear inequalities
\begin{eqnarray}
C_{ij}(\hat \delta) + \frac{\partial{C_{ij}(\hat \delta)}}{\partial{\delta_i}}(\delta_i - \hat \delta_i) \, + \, \frac{\partial{C_{ij}(\hat \delta)}}{\partial{\delta_j}} (\delta_i - \hat \delta_i) \quad \le \quad s_{ij}, && \quad \forall \ \hat \delta \in \R^n \label{cut}
\end{eqnarray}
\end{LE}

Constraints \eqref{cut} express the \textit{outer envelope} of the set described by \eqref{varf} \cite{BVbook}.  Any vector $\delta \in \cR^n$ which satisfies \eqref{varf} (for a given choice of $s_{ij}$) is guaranteed to satisfy \eqref{cut}.  Thus a a finite subset of the inequalities \eqref{cut}, used instead of \eqref{varf}, will give rise to a relaxation of the optimization problem and thus a lower bound on the optimal objective value.  Given Lemma \ref{outer} there are two ways to proceed, both motivated by the observation that at $\delta^*\in \R^n$, the most constraint inequality from among the set \eqref{cut} is that obtained by choosing $\hat \delta = \delta^*$.

First, one can use inequalities \eqref{cut} as \textit{cutting-planes} in the context of the ellipsoid method \cite{GLS}, obtaining a polynomial-time algorithm. A different way to proceed yields a numerically practicable algorithm.  For brevity, we will omit treatment of the generator conic constraints \eqref{gen1}, \eqref{gen2} (which are similarly handled). Denote by $F(\bar p, \alpha)$ the objective function in Eq.~\eqref{obj}.

\begin{center}
  \fbox{
    \begin{minipage}{0.95\linewidth}
      \begin{center}
      \begin{PRO}{CUTTING-PLANE ALGORITHM}\label{cutplanealg}\end{PRO}
      \end{center}
      {\bf Initialization:} The linear ``master'' system $A (\bar p, \alpha, \delta, \theta, s)^T \ge b$ is defined to include constraints \eqref{deltadef}-\eqref{sdef}.\\ \\
      {\bf Iterate:} \\ \\
      {\bf (1)} Solve $\min \{ F(\bar p, \alpha) \ : \ A(\bar p, \alpha, \delta, \theta, s)^T \ge b \}$. Let $(\bar p^*, \alpha^*, \delta^*, \theta^*, s^*)$ be an optimal solution.\\  \\
      {\bf (2)} If all conic constraints are satisfied up to numerical
tolerance by\\ $(\bar p^*, \alpha^*, \delta^*, \theta^*, s^*)$, {\bf Stop}.\\ \\
      {\bf (3)} If all chance constraints are satisfied up to numerical
tolerance by $(\bar p^*, \alpha^*)$, {\bf Stop}.\\ \\
      {\bf (4)} Otherwise, add to the master system the outer inequality \eqref{cut}
arising from that constraint (\ref{varf}) which is most violated.
\vspace{2pt}

\end{minipage}

  }
\vskip 3pt
\end{center}

As the algorithm iterates the master system represents a valid relaxation of the conic program \eqref{deltadef}-\eqref{varf}; thus the objective value of the solution computed in Step 1 is a valid lower bound on the value of problem.  Each problem solved in Step 1 is a linearly constrained, convex quadratic program.  Computational experiments involving large-scale realistic cases show that the algorithm is robust and rapidly converges to an optimum.

Note that Step 3 is not redundant.  The stopping condition in Step 2 may fail because the variance estimates are incorrect (too small), nevertheless the pair $(\bar p^*, \alpha^*)$ may already satisfy the chance constraints.  Checking that it does, for a given line $(i,j)$, is straightforward since the flow $\bmm{f_{ij}}$ is normally distributed (as noted in Lemma \ref{normalline}) and its mean and variance can be directly computed from $(\bar p^*, \alpha^*)$.

In our implementation termination is declared in Step 2 or Sep 3 when the corresponding constraint violation is less than $10^{-6}$. Table \ref{larger} displays typical performance of the cutting-plane algorithm on (comparatively more difficult) large problem instances.  In the Table, 'Polish1'-'Polish3' are the three Polish cases included in MATPOWER \cite{MATPOWER} (in Polish1 we increased loads by $30\%$). All Polish cases have uniform random costs on [0.5, 2.0] for each generator and ten arbitrarily chosen wind sources. The average wind power penetration for the four cases is 8.8\%, 3.0\%, 1.9\%, and 1.5\%.  'Iterations' is the number of linearly-constrained subproblems solved before the algorithm converges. 'Barrier iterations' is the total number of iterations of the barrier algorithm in CPLEX over all subproblems, and 'Time' is the total (wallclock) time required by the algorithm.  For each case, line tolerances are set to two standard deviations and generator tolerances three standard deviations.  These instances all prove unsolvable if directly tackled by CPLEX or Gurobi.

\begin{table}[h]
\centering
\caption{ Performance of cutting-plane method on a number of large cases.   }
\vskip 5pt
\begin{tabular}{l c c c c c c}
\hline
\hline
Case & Buses & Generators & Lines & Time (s) & Iterations & Barrier iterations \\
\hline
\hline
BPA & 2209 & 176 & 2866 & 5.51 & 2 & 256  \\
Polish1 & 2383 & 327 & 2896 & 13.64 & 13 & 535 \\
Polish2 & 2746 & 388 & 3514 & 30.16 & 25 & 1431 \\
Polish3 & 3120 & 349 & 3693  & 25.45 & 23 & 508 \\
\hline
\label{larger}
\end{tabular}
\end{table}

Table \ref{detailed} provides additional, typical numerical performance for the cutting-plane algorithm on an instance of the Polish grid model. Each row of Table \ref{detailed} shows that maximum relative error and objective value at the end of several iterations.  The total run-time was 25 seconds.  Note the ``flatness'' of the objective. This makes the problem nontrivial -- the challenge is to find a \textit{feasible} solution (with respect to the chance constraints); at the onset of the algorithm the computed solution is quite infeasible and it is this attribute that is quickly improved by the cutting-plane algorithm.

\begin{table}[h]
\centering
\caption{ Typical convergence behavior of cutting-plane algorithm on a large instance.}
\begin{tabular}{r c c}
\hline
Iteration & Max rel. error & Objective \\ \hline
1 & 1.2e-1 & 7.0933e6 \\
4 & 1.3e-3 & 7.0934e6 \\
7 & 1.9e-3 & 7.0934e6 \\
10 & 1.0e-4 & 7.0964e6 \\
12 & 8.9e-7 & 7.0965e6 \\ \hline
\label{detailed}
\end{tabular}
\end{table}

We note the (typical) small number of iterations needed to attain numerical convergence.  Thus at termination only a very small number of conic constraints \eqref{varf} have been incorporated into the master system.  This validates the expectation that only a small fraction of the conic constraints in CC-OPF are active at optimality.  The cutting-plane algorithm can be viewed as a procedure that opportunistically discovers these constraints.

\subsection{Data-robust chance constraints}
\label{subsec:ambiguous}

Above we developed a formulation for our chance-constrained OPF problem as
the conic program \eqref{obj}-\eqref{gen2}.  This approach assumed exact estimates for the mean $\mu_i$ and
the variance $\sigma^2_i$ of each wind source $\bm{\omega_i}$.  In practice
however the estimates at hand might be imprecise\footnote{In particular since they would effectively be computed in real time.}; consequently jeopardizing the usefulness of our conic program, henceforth termed the \textit{nominal} chance-constrained problem.  In particular, the performance of the control computed by the conic program might conceivably be sensitive to even small changes in the data.  We will deal with this issue in two complimentary ways.

\subsubsection{Out-of-sample analysis}

Our first approach is the out-of-sample tests, implemented experimentally in Section \ref{subsec:out-of-sample}. We assume that the $\mu_i$ and $\sigma^2_i$ have been mis-estimated and  explore the robustness of the affine control with respect to the estimation errors. The experiments of Section \ref{subsec:out-of-sample} show that the degradation of the chance constraints is small when small data errors are experienced. This empirical observation has a rigorous explanation discussed below.

Our chance constraints are represented by convex inequalities, for example in the case of $P(\bm{f_{ij}} > f_{ij}^{max}) < \epsilon_{ij}$ and $P(\bm{f_{ij}} < -f_{ij}^{max}) < \epsilon_{ij}$ we used e.g. \eqref{sdef} and \eqref{varf}.  Suppose that we have solved the chance-constrained problem assuming (incorrect) expectations $\mu_i$ and variances $\sigma_i^2$ ($i \in \cW$). Let $\tilde \mu_i$ and $ \tilde \sigma_i^2$ ($i \in \cW$) be the exact (or realized) values.  With some abuse of notation, we will write $\xi$ (resp., $\tilde \xi$) for the incorrect (exact) distribution. For a given line $(i,j)$, write:
\begin{eqnarray}
&&  m_{ij} \ = \ \mathbb{E_{\bm{ \xi}}} \bm{f_{ij}} \ = \
\beta_{ij}( [\breve \cB (\bar p  +  \mu - d) ]_j - [\breve \cB (\bar p  +  \mu - d) ]_j), \label{tildemean}\\
&&  \sigma^2_{ij} \ = \ var_{ \xi}\bm{f_{ij}} \ = \ \beta^2_{ij}\sum_{k \in \cW}  \sigma_k^2 (\pi_{i k} - \pi_{jk} - \delta_i + \delta_j)^2. \label{tildevar}
\end{eqnarray}
and similarly,
\begin{eqnarray}
&& \tilde m_{ij} \ = \ \mathbb{E_{\bm{\tilde \xi}}} \bm{f_{ij}} \ = \
\beta_{ij}( [\breve \cB (\bar p  + \tilde \mu - d) ]_j - [\breve \cB (\bar p  + \tilde \mu - d) ]_j), \label{tildemean2}\\
&& \tilde \sigma^2_{ij} \ = \ var_{\tilde \xi}\bm{f_{ij}} \ = \ \beta^2_{ij}\sum_{k \in \cW} \tilde \sigma_k^2 (\pi_{i k} - \pi_{jk} - \delta_i + \delta_j)^2. \label{tildevar2}
\end{eqnarray}
Using \eqref{tildemean2} and \eqref{tildevar2} we see that the ``true'' probability $P(\bm{f_{ij}} > f_{ij}^{max})$ is that value $\tilde \epsilon$ such that
\begin{eqnarray}
&& \tilde m_{ij} + \eta(\tilde \epsilon) \tilde\sigma_{ij} = u_{ij}. \label{tildepsilon}
\end{eqnarray}
We wish to evaluate how much \textit{larger} this (realized) value $\tilde \epsilon$ is compared with the target value $\epsilon_{ij}$ which was the goal in the chance-constrained computation. We will do this assuming that the estimation errors are small, more precisely, there exist  nonnegative constants $M$ and $V$ such that
\begin{eqnarray}
 &&  \ \ \forall i \in \cW, \quad |\tilde \mu_i- \mu_i| < M \mu_i \quad \ \mbox{and} \quad \ |\tilde \sigma^2_i - \sigma^2_i| < V^2\sigma^2_i.
\end{eqnarray}
Considering Eq.~\eqref{tildepsilon}, we see that for data errors of a given magnitude, $\tilde \epsilon$ is maximized when $\tilde m_{ij}$ and $\tilde \sigma_{ij}$
are maximized. Further, considering Eqs.~\eqref{tildemean2} and Eqs.~\eqref{tildevar2} we see that to first order $\tilde m_{ij} \le m_{ij} + O(M)$, and that   $\tilde \sigma^2_{ij} \le (1 + V^2) \sigma^2_{ij}$.  From these two observations and eq.~\eqref{tildepsilon} we obtain
\begin{eqnarray}
\tilde \epsilon & = & \frac{1}{\sqrt{2 \pi}} \int_{\eta(\tilde \epsilon)}^{+ \infty} e^{-x^2/2} dx \ = \ \epsilon_{ij} + \frac{1}{\sqrt{2 \pi}} \int_{\eta(\tilde \epsilon)}^{\eta(\epsilon_{ij})} e^{-x^2/2} dx \nonumber \\
&& = \epsilon_{ij} + \frac{1}{\sqrt{2 \pi}} e^{-\frac{\eta(\epsilon_{ij})^2}{2}}[ \eta(\epsilon_{ij}) - \eta(\tilde \epsilon) ] \ + \ \mbox{smaller order errors}. \label{firstestimate}
\end{eqnarray}
The quantity in brackets in the right-hand side of \eqref{firstestimate} equals
\begin{eqnarray}
&& \eta(\epsilon_{ij}) - \frac{u_{ij} -  m_{ij}}{\tilde \sigma_{ij}} \, + \, \frac{\tilde m_{ij} - m_{ij}}{\tilde \sigma_{ij}} \ < \\
&& \eta(\epsilon_{ij}) \left(1 - \frac{1}{\sqrt{1 + V^2}} \right) \, + \, \frac{\tilde m_{ij} - m_{ij}}{\tilde \sigma_{ij}} \ < \ \\ && \eta(\epsilon_{ij}) O(V)  \, + \, \frac{O(M)}{\sigma_{ij}}. \label{brack}
\end{eqnarray}
Using \eqref{firstestimate} and \eqref{brack} we obtain that
$\tilde \epsilon  = \epsilon_{ij} + O(V) + O(M)$, where the ``O'' notation
``hides'' solution-dependent constants.

\subsubsection{Efficiently solvable data-robust formulations}

As the preceding analysis makes clear, the constants $M$ and $V$ may need to be quite small, for example if the $\sigma_{ij}$ are small.  We thus seek a better guarantee of robustness. This justifies our second approach discussed below -- to develop a version of CC-OPF which is methodologically guaranteed to be insensitive to data errors. This approach puts CC-OPF within the scope of robust optimization; to be more precise we will be solving an ambiguous chance-constrained problem in the language of \cite{ambiguous}.

We will write each $\mu_i$ in the form $\bar \mu_i + r_i$, where the $\bar \mu_i$ are point estimates of the $\mu_i$ and the $r_i$ are ``errors'' which are constrained to lie in some fixed set $\cM$ with $0 \in \cM$.  Likewise, we assume that there is a set $\cS \subseteq \cR^{\cW}$ including $0$, such that each $\sigma_i^2$ is of the form $\bar \sigma_i^2 + v_i$ where the vector of errors $v_i$ belongs to $\cS$. As an example for how to construct  $\cM$ or $\cS$, one can use the following set parameterized by values $0 < \gamma_i$ and $0 < \Gamma$:
\begin{eqnarray}
U(\gamma, \Gamma) & = & \left\{ r \in \R^{\cW} \, : \, | r_i| \, \le \, \gamma_i \,\, \forall i \in \cW, \ \ \sum_{i \in \cW} \frac{| r_i|}{\gamma_i}  \le   \Gamma \right\}.
\end{eqnarray}
This set was introduced in \cite{bertsim}.  Another candidate is the ellipsoidal set
\begin{eqnarray}
E(A,b) & = & \left\{ r \in \R^{\cW} \, : \, \, r^T A r \le b \, \right\}, \label{ellipsoidal}
\end{eqnarray}
where $A \in \R^{\cW \times \cW}$ is positive-definite and $b \ge 0$ is a real; see \cite{nemben}, \cite{dongarud}. We can now formally proceed as follows: \\

\noindent {\bf Definitions.} Let the estimates $\bar \mu$, $\bar \sigma^2$ and the sets $\cM$ and $\cS$ be given.\\

\noindent {\bf 1.} Let $\cD = \cD( \bar \mu, \cM, \bar \sigma^2, \cS)$ be the set of vectors of random variables such that for each pair $\pi \in \cM$ and $v \in \cS$ there is an element $\bm{\xi} \in \cD$ with coordinate-wise mean $\bar \mu + \pi$ and variance $\bar \sigma^2 + s$.

\noindent {\bf 2.} A pair $\bar p, \alpha$ is called \textit{robust} with respect to the pair $\cM$, $\cS$, if for each line $(i,j)$
\begin{eqnarray}
&& \max_{\xi \in \cD} P(\bm{f_{ij}} > f^{max}_{ij}) \ < \ \epsilon_{ij}, \ \ \mbox{and} \label{robust1}\\
&& \max_{\xi \in \cD} P(\bm{f_{ij}} < -f^{max}_{ij}) \ < \ \epsilon_{ij}. \label{robust2}
\end{eqnarray}

Our task will be to replace, in our optimization problem formulation, the chance constraint \eqref{chance_constraint} with one asking for robustness as in Eqs.~\eqref{robust1}-\eqref{robust2}\footnote{And likewise with the generator chance-constraints \eqref{gen1}, \eqref{gen2}.}.  If the uncertainty sets $\cM$ and $\cS$ consist of a single point estimate each, then we recover the nominal chance-constrained problem we discussed above.  As $\cM$ or $\cS$ become larger, the robust approach becomes more insensitive to estimation errors, albeit at the cost of becoming more conservative. Thus, a reasonable balance should be attained by choosing $\cM$ and $\cS$ small but of positive measure, thereby preventing trivial sensitivity of the control to changes in the data.

To indicate our approach, we focus on one of the statements for our chance constraint for lines.  Using Lemma \ref{normalline}, the robustness criterion \eqref{robust1} applied to a line $(i,j)$ requires
\begin{eqnarray}
\beta_{ij} \max_{ \xi \in \cD}  \left \{ (\bar \theta_i - \bar \theta_j) + \eta(\epsilon_{ij}) \left[ \sum_{k \in \cW} \sigma_k^2 (\pi_{i k} - \pi_{jk} - \delta_i + \delta_j)^2\right]^{1/2} \right\}  & \le & f^{max}_{ij}. \label{ambig}
\end{eqnarray}
We can see that \eqref{ambig} consists of a (possibly infinite) set of
convex constraints; thus the data-robust chance-constrained problem is a
convex problem.  However we need to exploit this fact in a computationally favorable manner.  To motivate our approach, we will first show that a methodology based on traditional robust optimization techniques will (likely) not succeed at this task.  We will then describe our method.

Considering the expression in brackets in the left-hand side of Eq.~\eqref{ambig} we recall that $\bar \theta  = \bar \theta(\mu) =  \breve \cB (\bar p  + \mu - d)$ depends on the uncertainty set $\cM$, in fact we can write
\begin{eqnarray}
\bar \theta & =  & \breve \cB (\bar p  + \bar \mu - d) \ + \breve \cB r \ = \ \bar \theta(\bar \mu) \ + \breve \cB r,
\end{eqnarray}
where $r \in \cM$. Also, the $\pi$ are constants, and  $\delta = \breve \cB \alpha$ which is deterministic though dependent on our decision variables $\alpha$. In summary, constraint \eqref{ambig} can be rewritten as
\begin{eqnarray}
&& \beta_{ij} (\bar \theta(\bar \mu)_i - \bar \theta(\bar \mu)_j) \ + \ \beta_{ij} \max_{r \in \cM} \left\{ e^T_{ij}\breve \cB r \right\}\ + \ \nonumber \\\ && \eta(\epsilon_{ij}) \beta_{ij} \left[ \sum_{k \in \cW} \bar \sigma_k^2 (\pi_{i k} - \pi_{jk} - \delta_i + \delta_j)^2 \ + \ \max_{v \in \cS} \left\{ \sum_{k \in \cW} v_k (\pi_{i k} - \pi_{jk} - \delta_i + \delta_j)^2  \right\} \right]^{1/2}\nonumber \\   &&  \le \ f^{max}_{ij}. \label{ambig2}
\end{eqnarray}
where $e_{ij}^T \in \R^n$ is the vector with a +1 entry in position $i$, a -1 entry in position $j$ and 0 elsewhere.   Note that if in the left-hand side of \eqref{ambig2} we ignore second term and the second term inside the square brackets,  we obtain the nominal (i.e. non-robust) version of chance-constraint \eqref{robust1}; see Eqs.~\eqref{sdef}, \eqref{varf}. A similar constraint is obtained from Eqs.~\eqref{robust2}.

Constraint \eqref{ambig2} can be incorporated into a formulation provided we can appropriately represent the two maxima. Here, note that
\begin{eqnarray}
&& m_{ij} \ \doteq \ \max_{r \in \cM} \left\{ e^T_{ij}\breve \cB r \right\}, \label{mijdef}
\end{eqnarray}
is independent of all variables and can be solved beforehand for all lines $(i,j)$;  when $\cM$ is of the form $U(\gamma, \Gamma)$ given above this is a linear programming problem and when $\cM = E(A, b)$ the task amounts to finding a point in the boundary of an ellipsoid with normal parallel to a given vector and thus requires solving a linear system of equations.  We will likewise define a quantity $m_{ji}$ using $e_{ji}$ instead of $e_{ij}$.

It is the second maximization that presents some challenges.
\begin{LE}\label{lindual} Suppose $\cS = U(\gamma, \Gamma)$, and suppose a vector $\delta \in R^n$ is given.  Then
\begin{eqnarray}
\max_{v \in \cS} \left\{ \sum_{k \in \cW} v_k (\pi_{i k} - \pi_{jk} - \delta_i + \delta_j)^2  \right\} \quad &= & \ \min \ \ \Gamma a + \sum_{k \in \cW} \gamma_k b_k \nonumber \\ \  \mbox{s.t.}&& \ a + b_k \, \ge \, (\pi_{i k} - \pi_{jk} - \delta_i + \delta_j)^2 \ \ \forall \, k \in \cW, \nonumber \\
 && \ \ b_k \ge 0 \ \ \forall \, k \in \cW; \ \ a \ge 0. \nonumber
\end{eqnarray}
\end{LE}
\noindent {\em Proof sketch.}
Since $(\pi_{i k} - \pi_{jk} - \delta_i + \delta_j)^2 \ge 0$ without loss of generality in the maximum we will have that $v_k \ge 0$ for all $k$.  Strong linear programming duality now gives the result. \QED \\

The use of linear programming duality as in Lemma \ref{lindual} is key in the context of sets of the form $U(\gamma, \Gamma)$. In the case of an ellipsoidal set $E(A,b)$ as in \eqref{ellipsoidal} there is an analogue to Lemma \ref{lindual} that instead uses the S-Lemma \cite{BVbook}, \cite{dongarud}.  In the standard robust optimization approach, Lemma \ref{lindual} would be leveraged to produce a result of the following type:
\begin{LE}\label{formlemma} Suppose $\cS = U(\gamma, \Gamma)$.  The data-robust chance-constrained problem is obtained by replacing for each line $(i,j)$ constraints \eqref{sdef}, \eqref{varf} of the nominal formulation with
\begin{eqnarray}
&& \beta_{ij} (\bar \theta(\bar \mu)_i - \bar \theta(\bar \mu)_j) \ + \ \beta_{ij} m_{ij} \ + \ \beta_{ij} \eta(\epsilon_{ij}) \, s_{ij} \ \leq \ f^{max}_{ij}  \label{u1} \\
&& \beta_{ij} (\bar \theta(\bar \mu)_j - \bar \theta(\bar \mu)_i) \ + \ \beta_{ij} m_{ji} \ + \ \beta_{ij} \eta(\epsilon_{ij}) \, s_{ij} \ \leq \ f^{max}_{ij}  \label{u2} \\
&&  \left[\sum_{k \in \cW} \bar \sigma_k^2 (\pi_{i k} - \pi_{jk} - \delta_i + \delta_j)^2 \ + \ \Gamma a^{\{i,j\}} \ + \ \sum_{k \in \cW} b^{\{i,j\}}_k \right]^{1/2} \ \le \ s_{ij} \label{varf2ijk}\\
&& (\pi_{i k} - \pi_{jk} - \delta_i + \delta_j)^2 \ - \ a^{\{i,j\}} - b^{\{i,j\}}_k \ \le \ 0 \ \ \forall k \in \cW  \label{varf3ijk}\\
&& \ \ b^{\{i,j\}}_k \ge 0 \ \ \forall k \in \cW; \ \ a^{\{i,j\}} \ge 0.
\end{eqnarray}
Here, $s_{i,j}, \, a^{\{i,j\}}$ and $b^{\{i,j\}}_k$ ($k \in \cW$) are additional variables.
\end{LE}
\textit{Proof sketch.} Without loss of generality at the optimum the $a^{\{i,j\}}$ and $b_k^{\{i,j\}}$ are chosen so as to minimize the left-hand side of \eqref{varf2ijk} subject to all other variables held fixed, thereby obtaining the ``min'' in Lemma
\ref{lindual}. \QED

Lemma \ref{formlemma} points out the difficulty that the standard robust optimization approach would engender in the context of our problem.  First, the number of constraints \eqref{varf2ijk} and \eqref{varf3ijk} is \textit{large}: it equals $|\cE|(1 + |\cW|)$ and thus in the case of a large transmission system it could approach many tens of thousands (or more).  Thus, even though we obtain a \textit{compact} formulation (i.e., of polynomial size) it is likely to be proven too large for present-day solvers. But there is a  second and more fundamental problem: constraint \eqref{varf2ijk} \textit{is not convex}.  This is a significant methodological difficulty.    A similar set of hurdles arises when using uncertainty sets $E(A,b)$.

\subsubsection{Efficient solution of the data-robust problem using cutting planes}
To derive an algorithm for the data-robust chance-constrained problem that
(a) has some theoretical justification and (b) can prove numerically tractable,
we return to inequality \eqref{ambig2}, and note that if we replace the set
$\cS$ with a finite subset $\tilde \cS \subseteq \cS$ we obtain a valid relaxation.
In other words
\begin{eqnarray}
&& \beta_{ij} (\bar \theta(\bar \mu)_i - \bar \theta(\bar \mu)_j) \ + \ \beta_{ij} m_{ij} \ + \ \beta_{ij} \eta(\epsilon_{ij}) \, s_{ij} \ \leq \ f^{max}_{ij}  \label{u1prime} \\
&& \beta_{ij} (\bar \theta(\bar \mu)_j - \bar \theta(\bar \mu)_i) \ + \ \beta_{ij} m_{ji} \ + \ \beta_{ij} \eta(\epsilon_{ij}) \, s_{ij} \ \leq \ f^{max}_{ij}  \label{u2prime} \\
&&  \max_{v \in \tilde \cS} \left[  \sum_{k \in \cW} (\bar \sigma_k^2 + v_k) (\pi_{i k} - \pi_{jk} - \delta_i + \delta_j)^2  \right]^{1/2} \ \le \ s_{ij} \label{finite}
\end{eqnarray}
constitutes a valid relaxation of the chance constraints \eqref{sdef}, \eqref{varf} of the nominal formulation for each line $(i,j)$ for any given
finite $\tilde \cS \subseteq \cS$.    This observation can be used to
formally obtain a polynomial-time algorithm for the data-robust problem in
the cases of interest. For a given $v \in \tilde \cS$ let
$V_{i,j}^v(\delta)$ denote the expression inside the ``max'' in \eqref{finite}.  For
completeness, we state the following:
\begin{LE} \label{ellipsoid} In the case of uncertainty sets of the
form $U(\gamma, \Gamma)$ or $E(A,b)$ the data-robust chance-constrained
problem can be solved in polynomial time. \end{LE} \\
\noindent {\em Proof.}  Suppose we are given quantities $\hat \delta_i$ for
each $i \in \cV$ and $\hat s_{ij}$ for each line $(i,j)$.  Then as argued
before, if $\cS$ is of the form $U(\gamma, \Gamma)$ or $E(A,b)$ one can check
in polynomial time whether $\max_{v \in \cS}\left\{ V_{i,j}^v(\delta) \right\} \le  \hat s_{ij}$.  If the condition is violated for $\bar v \in \cS$ then
trivially
\begin{eqnarray}
L^{\bar v}_{ij}(\hat \delta) + \frac{\partial{L^{\bar v}_{ij}(\hat \delta)}}{\partial{\delta_i}}(\delta_i - \hat \delta_i) \, + \, \frac{\partial{L^{\bar v}_{ij}(\hat \delta)}}{\partial{\delta_j}} (\delta_i - \hat \delta_i) \quad \le \quad s_{ij}, &&  \label{viol2}
\end{eqnarray}
is violated at $\hat \delta, \hat s$.  Since \eqref{viol2} is valid
for the data-robust chance-constrained problem (by convexity) the result follows
by relying on the ellipsoid method.  \QED \\

Lemma \ref{ellipsoid} describes a formally ``good'' algorithm.  For
computational purposes we instead would rely on a cutting-plane
algorithm much like Algorithm \ref{cutplanealg}, developed in Section \ref{subsec:cuttingplane} so as to handle the nominal chance-constrained problem.
In the data-robust setting the algorithm is almost identical,
except that given quantities $\hat \delta_i$ for
each $i \in \cV$ and $\hat s_{ij}$ for each line $(i,j)$ we proceed as in
the proof of Lemma \ref{ellipsoid} by finding a $v \in \cS$ that proves
most stringent for the pair $\hat \delta, \hat s$ and then add to the master
formulation constraint \eqref{viol2}. Details are left to the reader.

\section{Experiments/Results}
\label{sec:experiments}

Here we will describe qualitative
aspects of our affine control on small systems; in particular we focus
on the contrast between standard OPF and nominal CC-OPF, on problematic features
that can arise because of fluctuating wind sources and on out-of-sample
testing of the CC-OPF solution, including the analysis of non-Gaussian
distributions.  Some of our tests involve the BPA grid and Polish Grid, which are large;
we present additional set of tests to address the scalability of our solution methodology to the large cases.

\subsection{Failure of standard OPF}
Above (see Eq.~\eqref{OPF-opt0}) we introduced the so-called standard OPF method for setting traditional generator output levels.  When renewables are present, the natural extension of this approach would make use of some fixed estimate of output (e.g., mean output) and to handle fluctuations in renewable output through the same method used to deal with changes in load: ramping output of traditional generators  up or down in proportion to the  net increase or decrease in renewable output.   This feature could seamlessly be handled using today's control structure, with each generator's output adjusted at a fixed (preset) rate. For the sake of simplicity, in the experiments below we assume that all ramping rates are equal.

Different assumptions on these fixed rates will likely produce different numerical results; however, this general approach entails an inherent weakness. The key point here is that mean generator output levels \textit{as well as} in particular the ramping rates would be chosen \textit{without} considering the stochastic nature of the renewable output levels. Our experiments are designed to highlight the limitations of this ``risk-unaware'' approach.  In contrast, our CC-OPF produces control parameters (the $\bar p$ and the $\alpha$) that are risk-aware and, implicitly, also topology-aware -- in the sense of network proximity to wind farms.

\begin{figure} \centering
\includegraphics[width=0.9\textwidth]{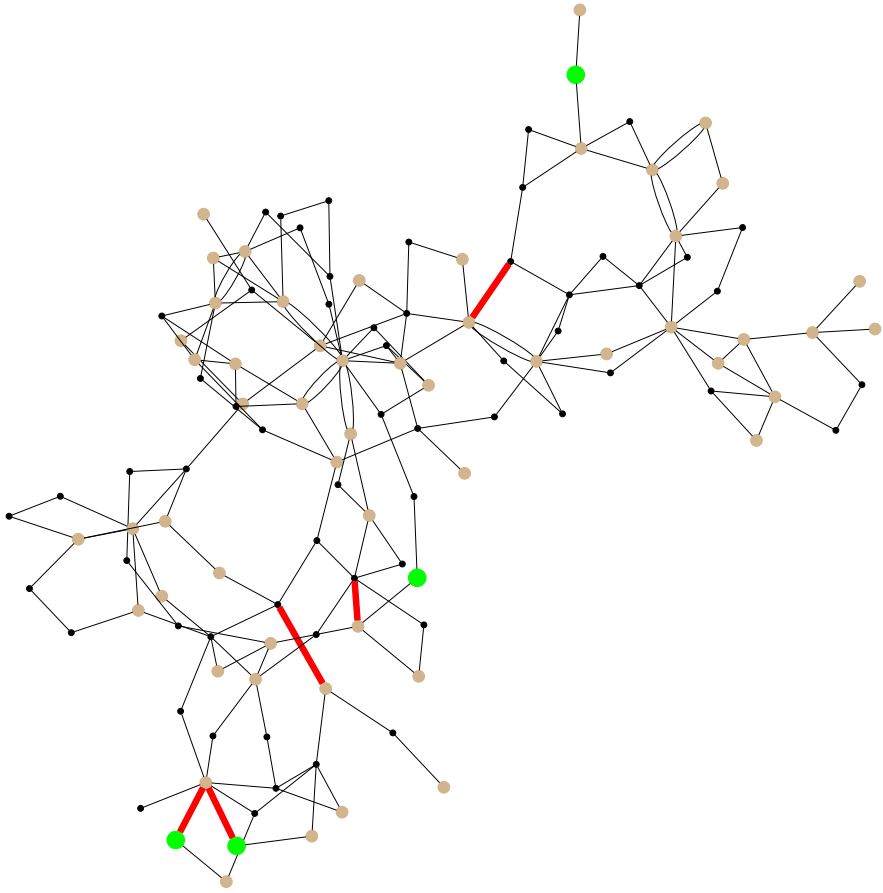} \caption{118-bus case
with four wind farms (green dots; brown are generators, black are loads). Shown
is the standard OPF solution against the average wind case with penetration of
5\%. Standard deviations of the wind are set to 30\% of the respective average
cases. Lines in red exceed their limit 8\% or more of the time.
\label{standard_bad}} \end{figure}

\begin{figure} \centering
\includegraphics[width=0.45\textwidth]{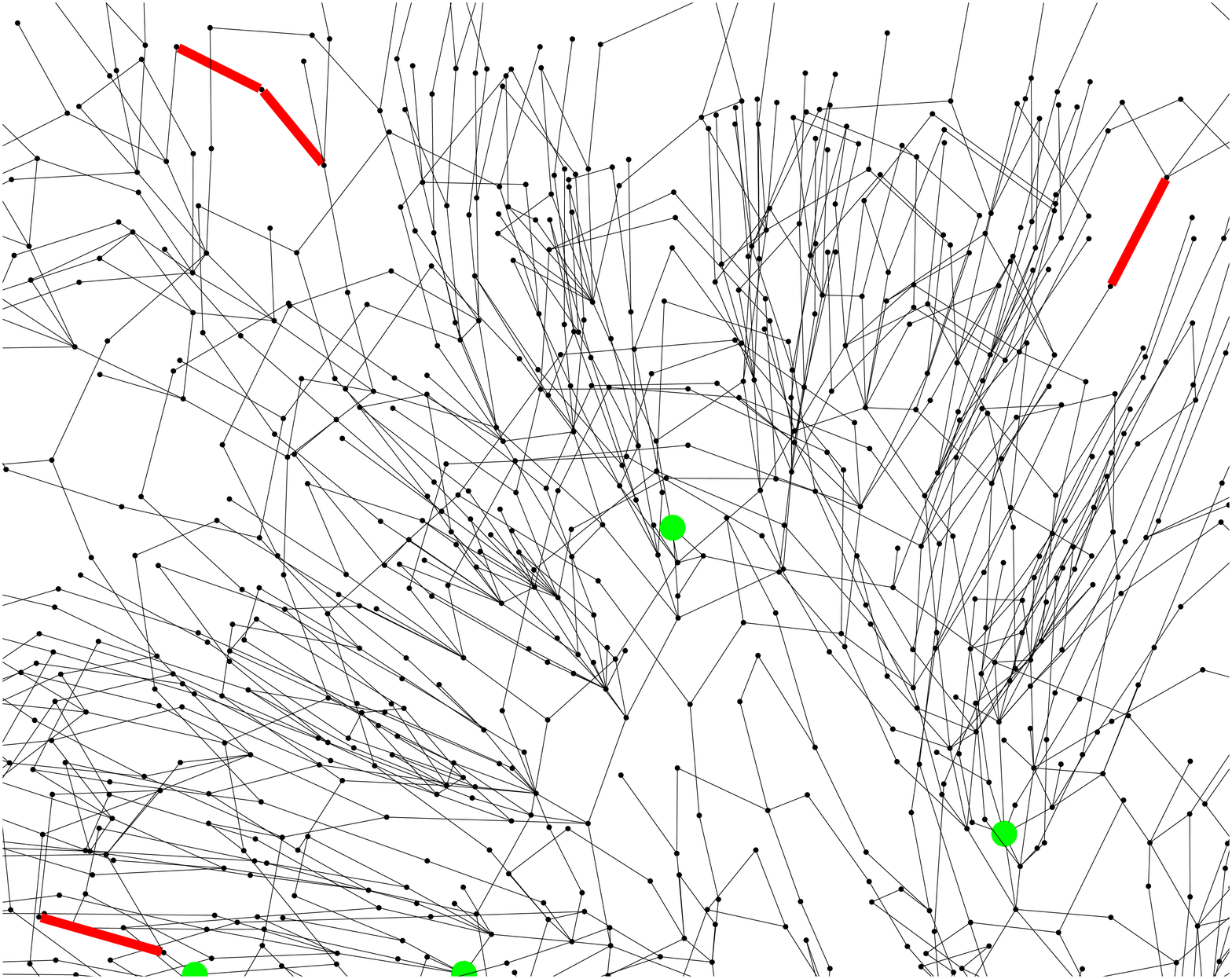}
\includegraphics[width=0.45\textwidth]{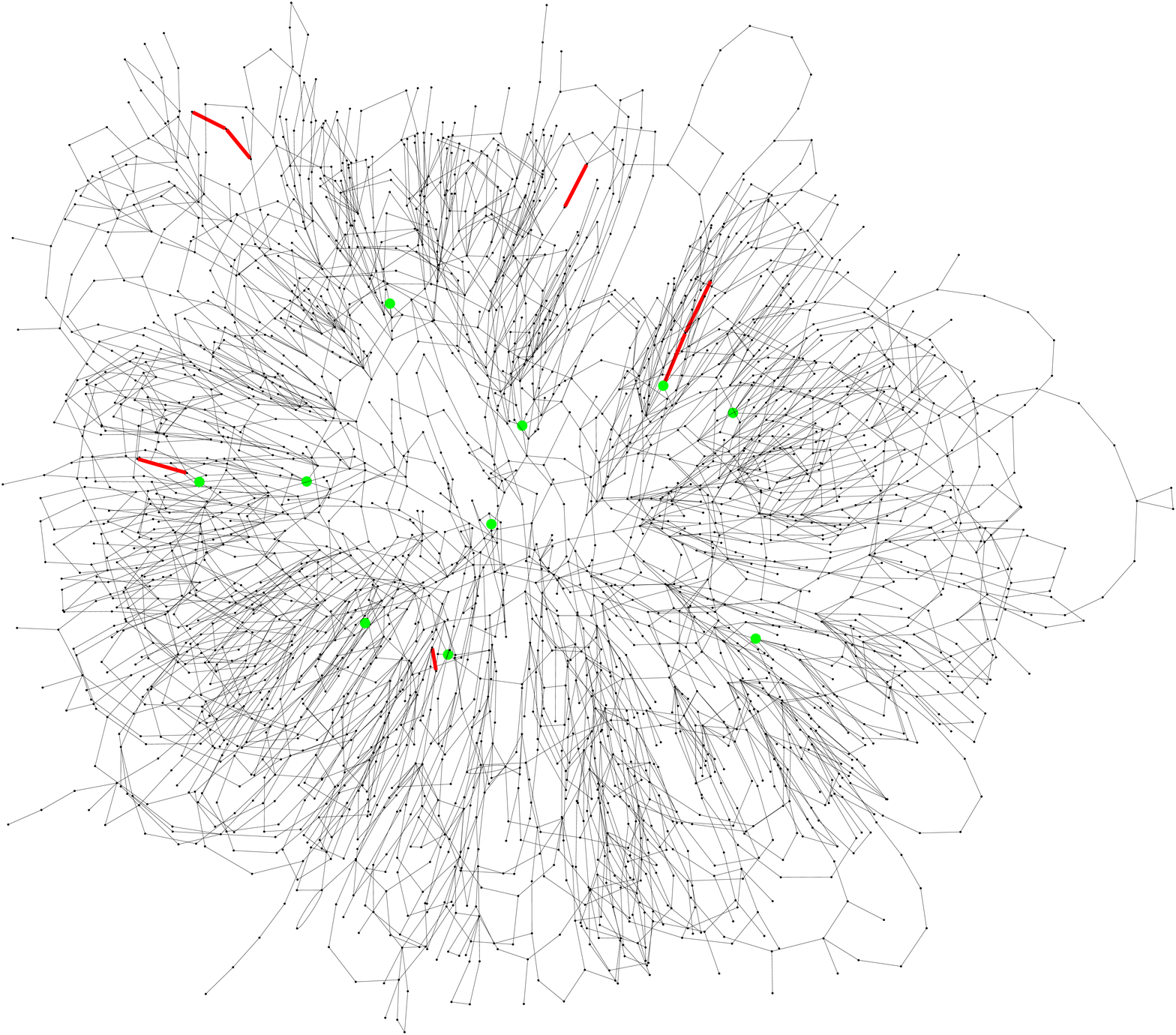}
\caption{Failure of the standard OPF shown for partial (left) and full (right)
snapshot of the standard OPF
solution on Polish grid from MATPOWER \cite{MATPOWER}. Color coding and conditions of the experiment are equivalent to these of
Fig.~\ref{standard_bad}.
\label{polish_bad}} \end{figure}

We first consider the IEEE 118-bus model with a quadratic cost function, and four sources of wind power added at arbitrary buses to meet 5\% of demand in the case of average wind.  The standard OPF solution is safely within the thermal capacity limits for all lines in the system.  Then we account for fluctuations in wind assuming Gaussian and site-independent fluctuations with standard deviations set to 30\% of the respective means. The results, which are shown in Fig.~\ref{standard_bad}, illustrate that under standard OPF five lines (marked in red) frequently become overloaded, exceeding their limits 8\% or more of the time. This situation translates into an unacceptably high risk of failure for any of the five red lines. This problem occurs for grids of all sizes; in Fig.~\ref{polish_bad} we show similar results on a 2746-bus Polish grid. In this case, after scaling up all loads by 10\% to simulate a more highly stressed system, we added wind power to ten buses for a total of 2\% penetration. The standard solution results in six lines exceeding their limits over 45\% of the time, and in one line over 10\% of the time.   For an additional and similar experiment using the Polish grid see Section \ref{subsec:scalability}.

\subsection{Cost of reliability under high wind penetration} 

The New York Times article ``Wind Energy Bumps Into Power Grid's Limits'' \cite{NYT2008} discusses how transmission line congestion has forced temporary shutdown of wind farms even during times of high wind. Our methodology suggests, as an alternative solution to curtailment of wind power, an appropriate reconfiguration of standard generators. If successful, this solution can use the available wind power without curtailment, and thus result in cheaper operating costs.

As a (crude) proxy for curtailment, we perform the following experiment, which considers different levels of renewable penetration. Here, the mean power outputs of the wind sources are kept in a fixed proportion to one another and proportionally scaled so as to vary total amount of penetration, and likewise with the standard deviations. First, we run our CC-OPF under a high penetration level.  Second, we add a 10\% buffer to the line limits and reduce wind penetration (i.e., curtail) so that under the \textit{standard} OPF  solution line overloads are reduced to an acceptable level. Assuming zero cost for wind power, the difference in cost for the high-penetration CC-OPF solution and the low-penetration standard solution are the savings produced by our model (generously, given the buffers).

\begin{figure} \centering
\includegraphics[width=0.45\textwidth]{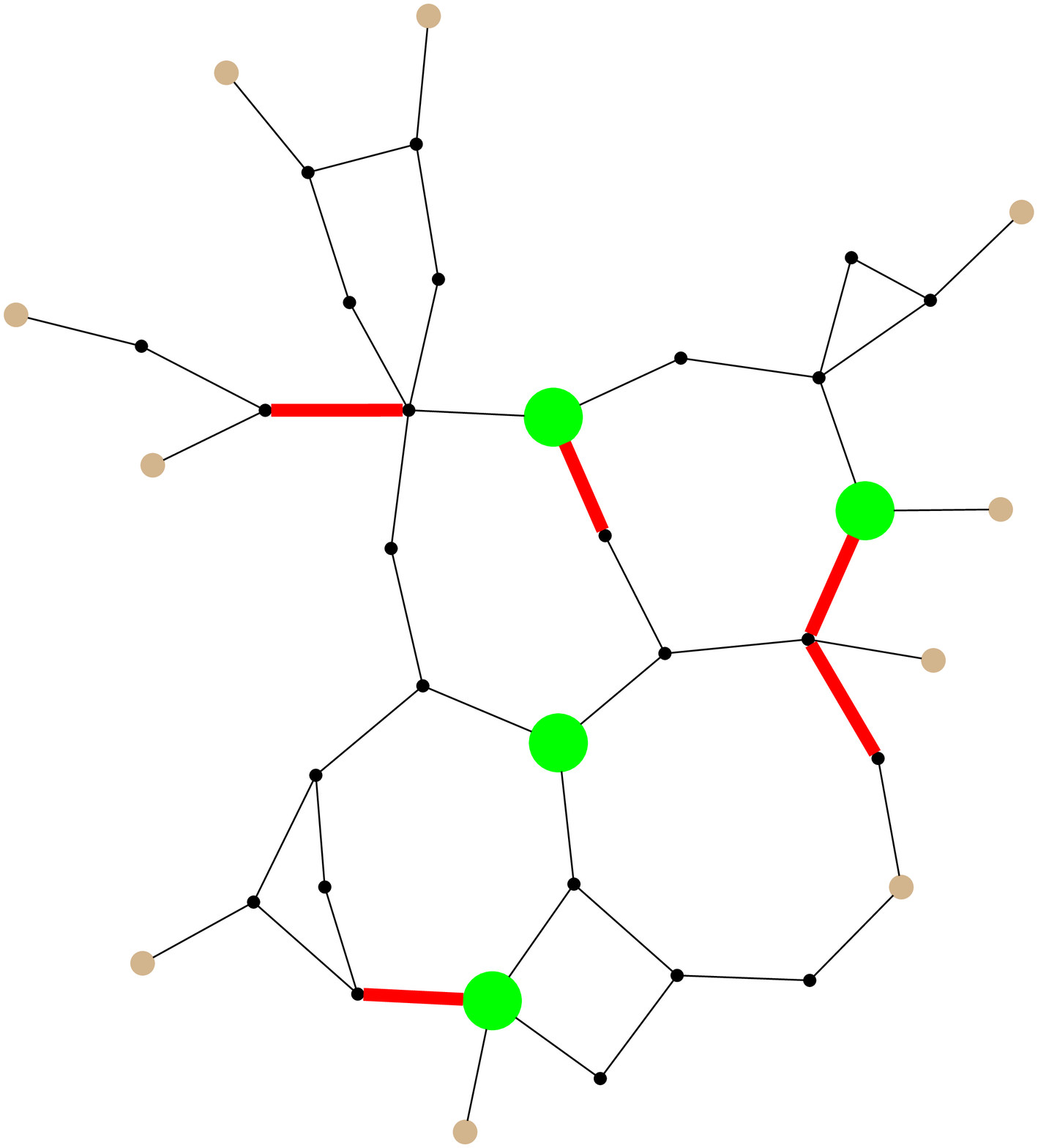}
\includegraphics[width=0.45\textwidth]{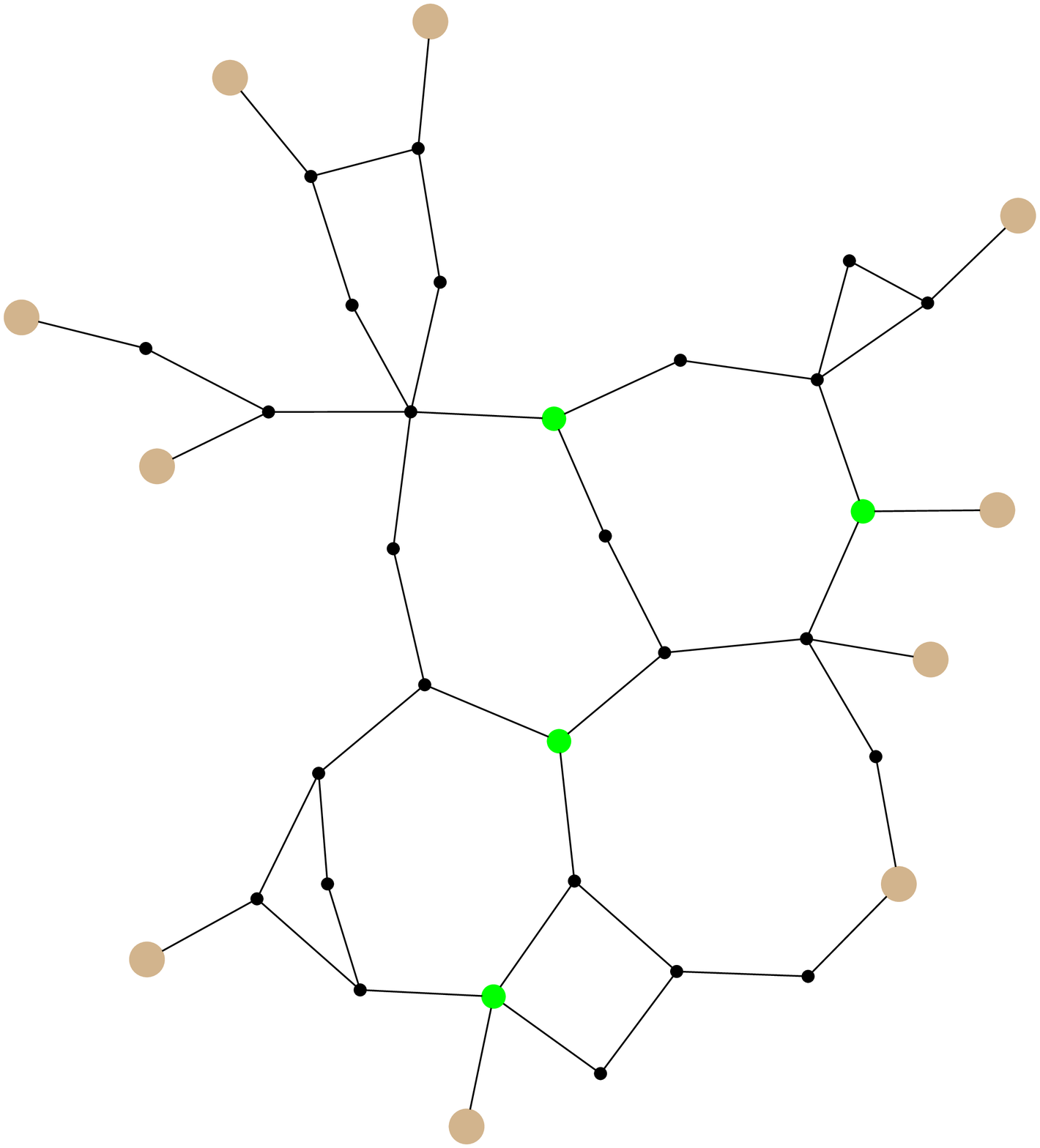} \caption{39-bus case
under standard solution. Even with a 10\% buffer on the line flow limits, five
lines exceed their limit over 5\% of the time with 30\% penetration (left). The
penetration must be decreased to 5\% before the lines are relieved, but at
great cost (right). The CC-OPF model is feasible for 30\% penetration at a cost
of 264,000. The standard solution at 5\% penetration costs 1,275,020 -- almost
5 times as much.  \label{cost}} \end{figure}

For the 39-bus case, our CC-OPF solution is feasible under 30\% of wind penetration, but the standard solution has 5 lines with excessive overloads, even when solved with the 10\% buffer. Reducing the penetration to 5\% relieves the lines, but more than quadruples(!) the cost over the CC-OPF solution. See Figure~\ref{cost}.  Note that this approach does not only show the advantage of the CC-OPF over standard OPF but also provides a quantitative measure of the advantage, thus placing a well-defined price tag on reliability.

\subsection{Non-locality} 

\begin{figure} \centering
\includegraphics[width=0.9\textwidth]{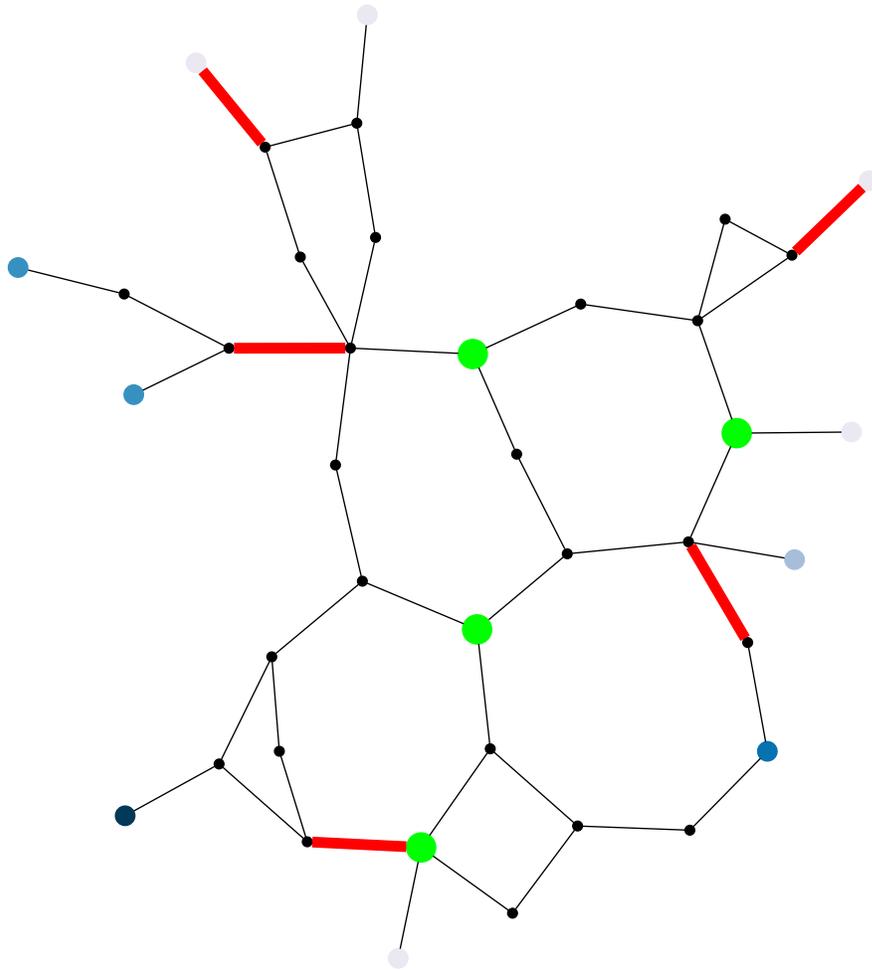} \caption{39-bus
case. Red lines indicate high probability of flow exceeding the limit under the
standard OPF solution. Generators are shades of blue, with darker shades
indicating greater absolute difference between the chance-constrained solution
and the standard solution.  \label{nonlocal}} \end{figure}

We have established that under fluctuating power generation, some lines may exceed their flow limits an unacceptable fraction of the time. Is there a simple solution to this problem, for instance, by carefully adjusting (a posteriori of the standard OPF) the outputs of the generators near the violated lines? The answer is no. Power systems exhibit significant non-local behavior. Consider Fig.~\ref{nonlocal}. In this example, the major differences in generator outputs between the standard OPF solution and our CC-OPF model's solution are not obviously associated with the different line violations. In general, it seems that it would be difficult to by-pass CC-OPF and make small local adjustments to relieve the stressed lines. On the positive side, even though CC-OPF is not local and requires a centralized computation,  it is also only slightly more difficult than the standard OPF in terms of implementation.

\subsection{Increasing penetration} 

Current planning for the power system in the United States calls for 30\% of wind energy penetration by 2030 \cite{20-2030}. Investments necessary to achieve this ambitious target may focus on both software (improving operations) and hardware (building new lines, sub-stations, etc), with the former obviously representing a much less expensive and thus economically attractive option.  Our CC-OPF solution contributes toward this option. A natural question that arises concerns the maximum level of penetration one can safely achieve by upgrading from the standard OPF to our CC-OPF.

\begin{figure} \centering
\includegraphics[width=0.3\textwidth]{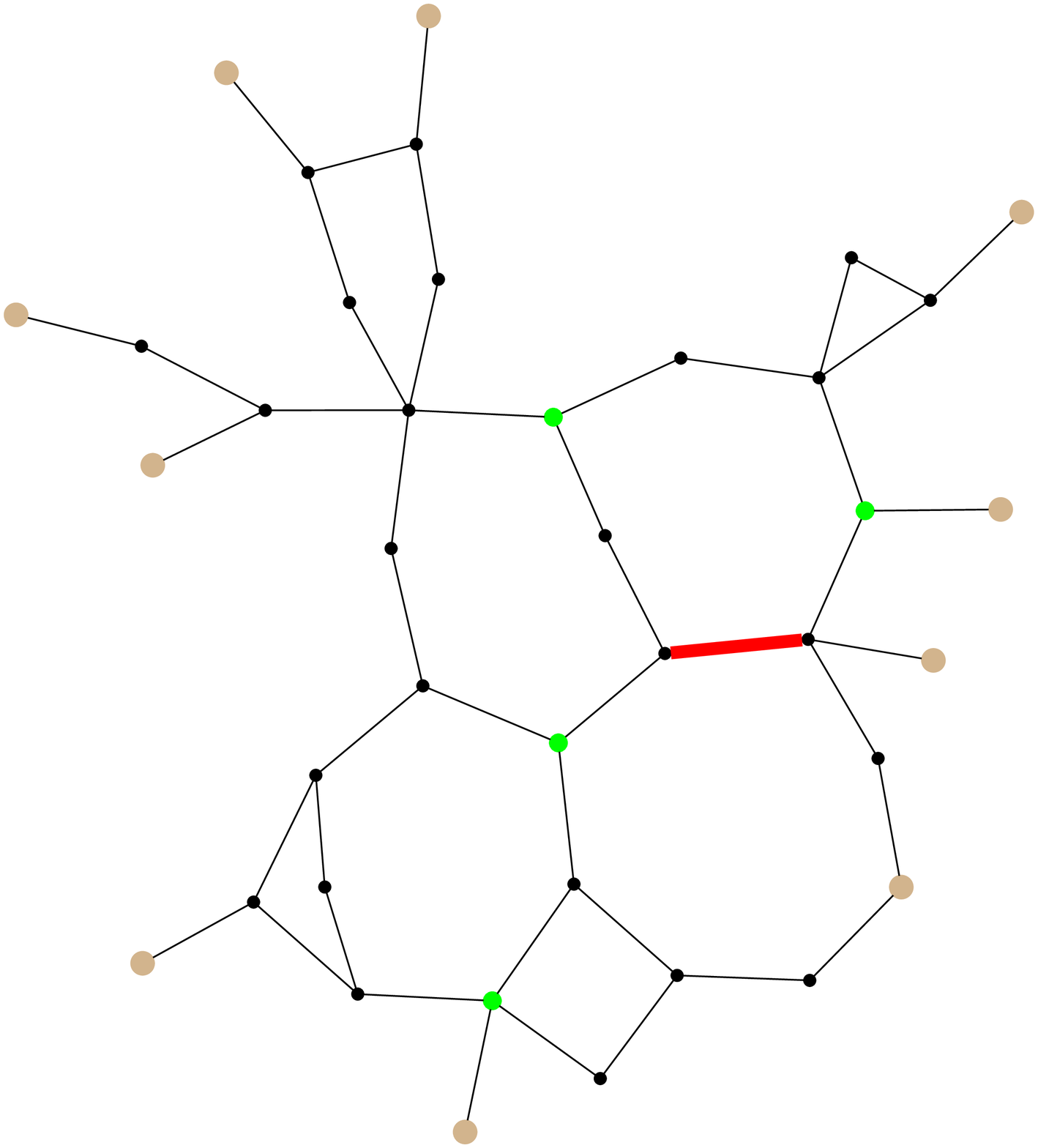}
\includegraphics[width=0.3\textwidth]{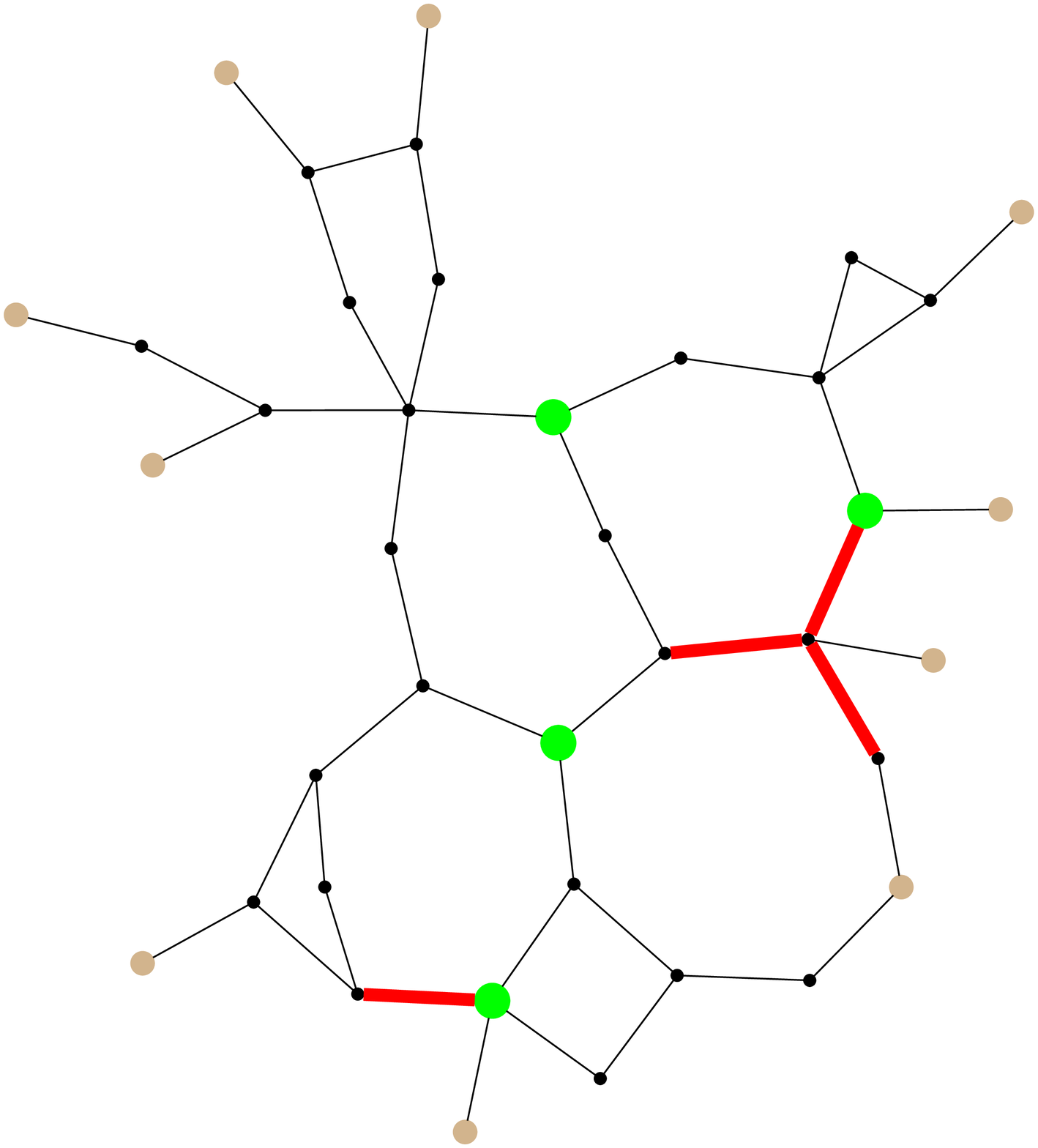}
\includegraphics[width=0.3\textwidth]{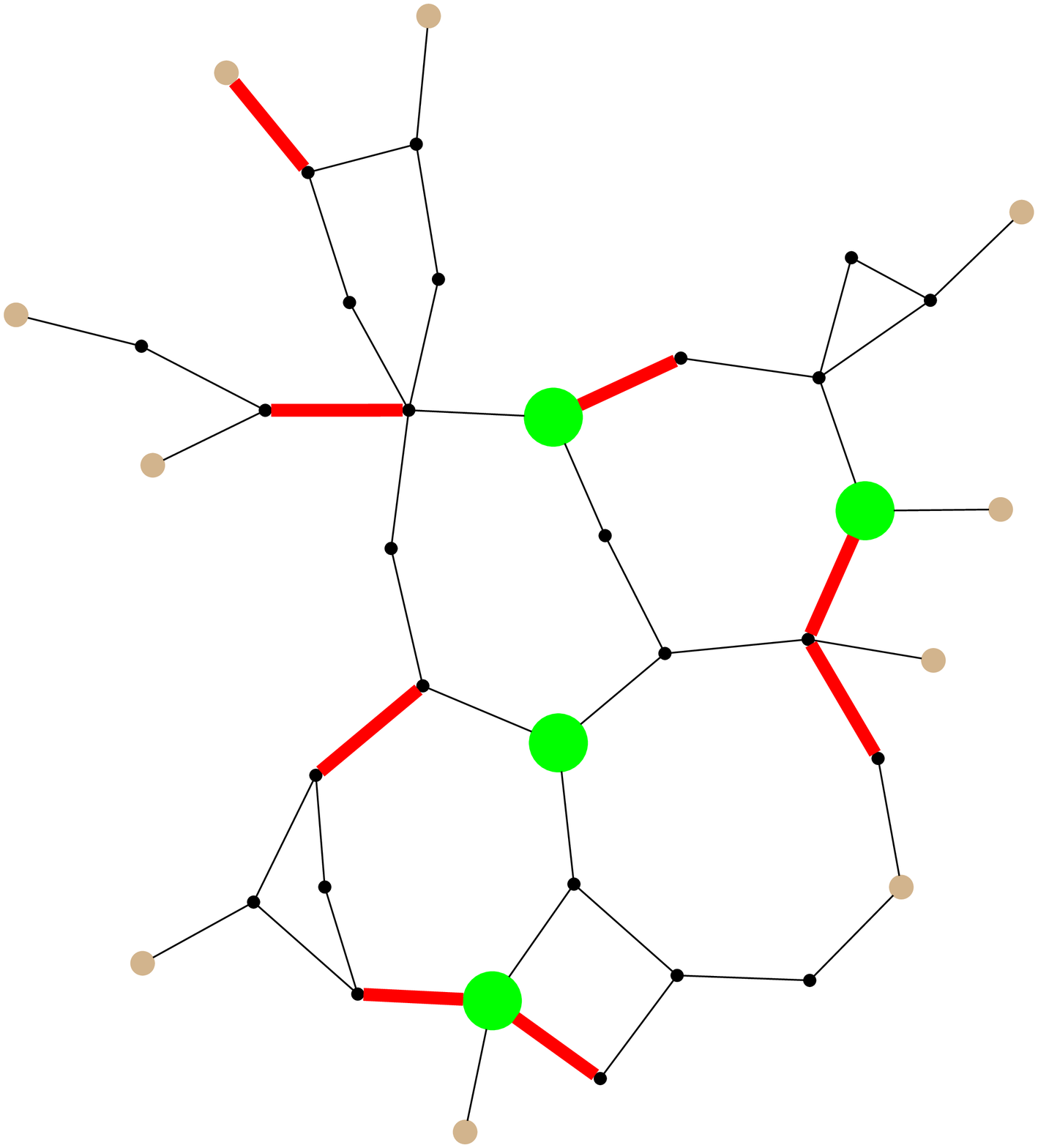} \caption{39-bus
case with four wind farms (green dots; brown are generators, black are loads).
Lines in red are at the maximum of $\epsilon = .02$ chance of exceeding their
limit. The three cases shown are left to right .1\%, 8\%, and 30\% average wind
penetration. With penetration beyond 30\% the problem becomes
infeasible.\label{penetration}} \end{figure}

To answer the question we consider the 39-bus New England system (from \cite{MATPOWER}) case with four wind generators added, and line flow limits scaled by .7 to simulate a heavily loaded system. The quadratic cost terms are set to rand(0,1) + .5. We fix the four wind generator average outputs in a ratio of 5/6/7/8 and standard deviations at 30\% of the mean. We first solve our model using $\epsilon = .02$ for each line and assuming zero wind power, and then increase total wind output until the optimization problem becomes infeasible. See Figure~\ref{penetration}.  This experiment illustrates that at least for the model considered, the 30\% of wind penetration with rather strict probabilistic guarantees enforced by our CC-OPF may be feasible, but in fact lies rather close to the dangerous threshold. To push penetration beyond the threshold is impossible without upgrading lines and investing in other (not related to wind farms themselves) hardware.

\subsection{Changing locations for wind farms} 

\begin{figure} \centering
\includegraphics[width=0.45\textwidth]{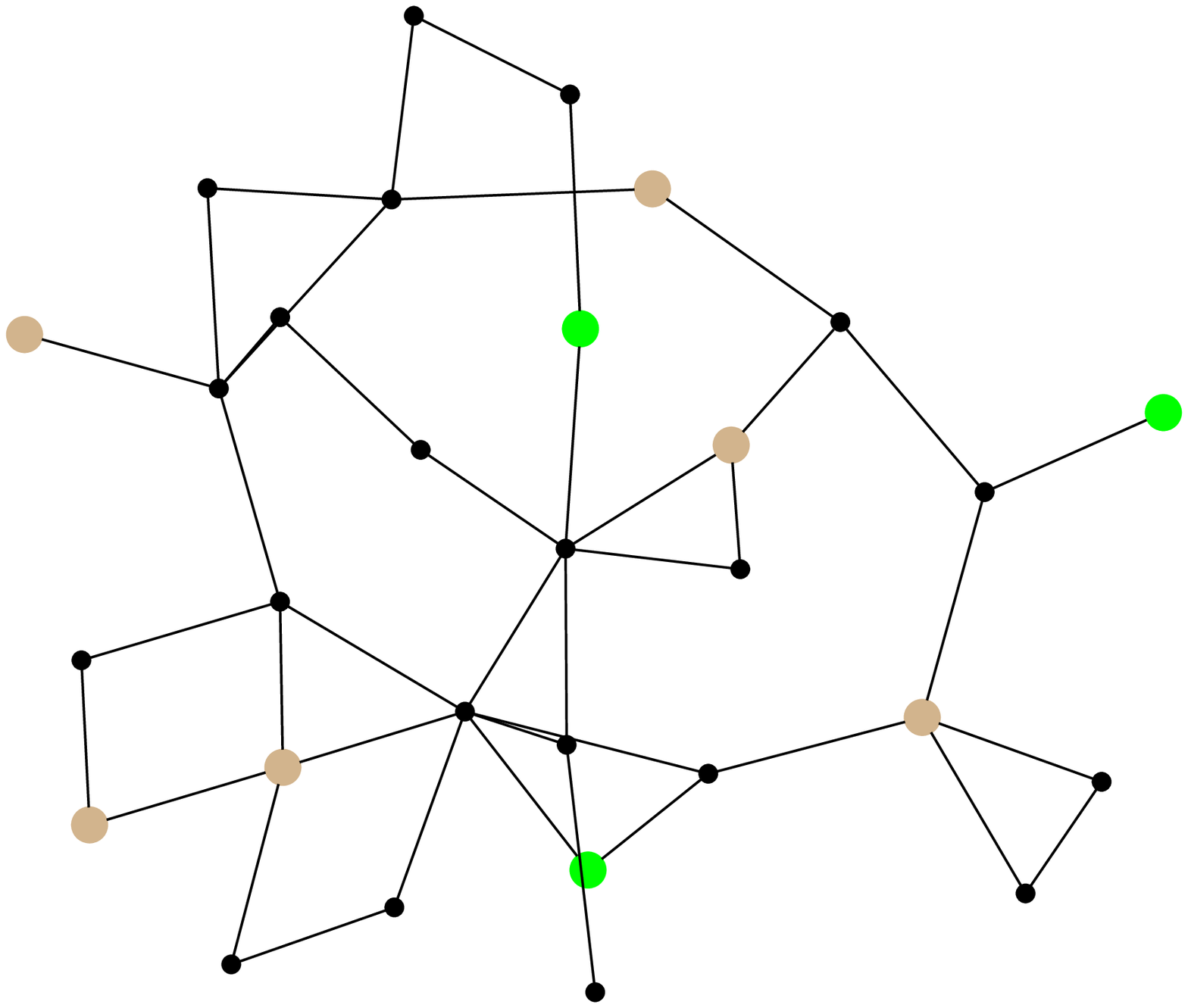}
\includegraphics[width=0.45\textwidth]{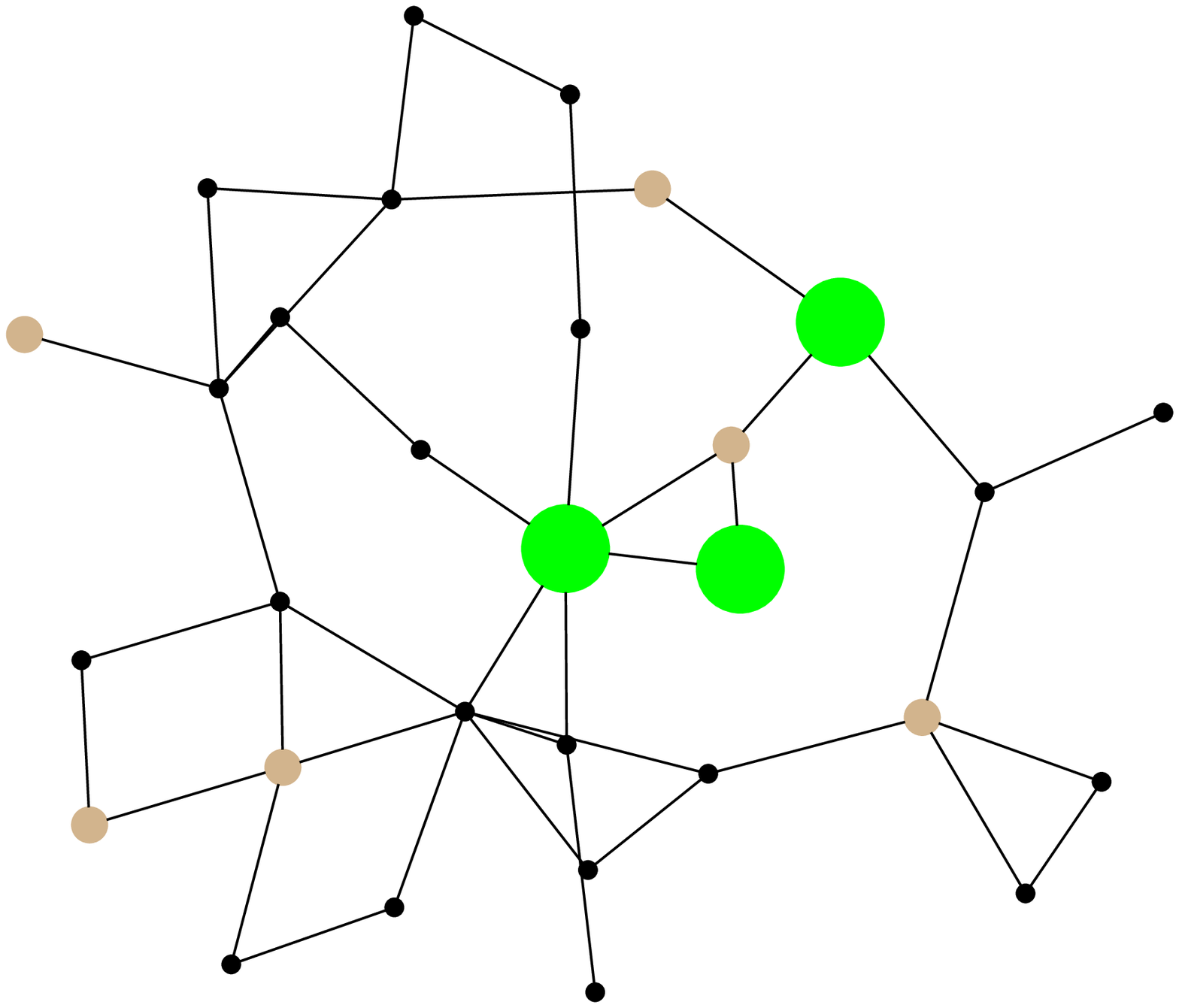} \caption{30-bus
case with three wind farms. The case on the left supports only up to 10\%
penetration before becoming infeasible, while the one on the right is feasible
for up to 55\% penetration.\label{diff_locations}} \end{figure}

In this example we consider the effect of changing locations of wind farms. We take the MATPOWER 30-bus case with line capacities scaled by .75 and add three wind farms with average power output in a ratio of 2/3/4 and standard deviations at 30\% of the average. Two choices of locations are shown in figure~\ref{diff_locations}. The first remains feasible for penetration up to 10\% while the second can withstand up to 55\% penetration. This experiment shows that choosing location of the wind farms is critical for achieving the ambitious goal of high renewable penetration.

\subsection{Reversal of line flows} 

\begin{figure} \centering
\includegraphics[width=0.45\textwidth]{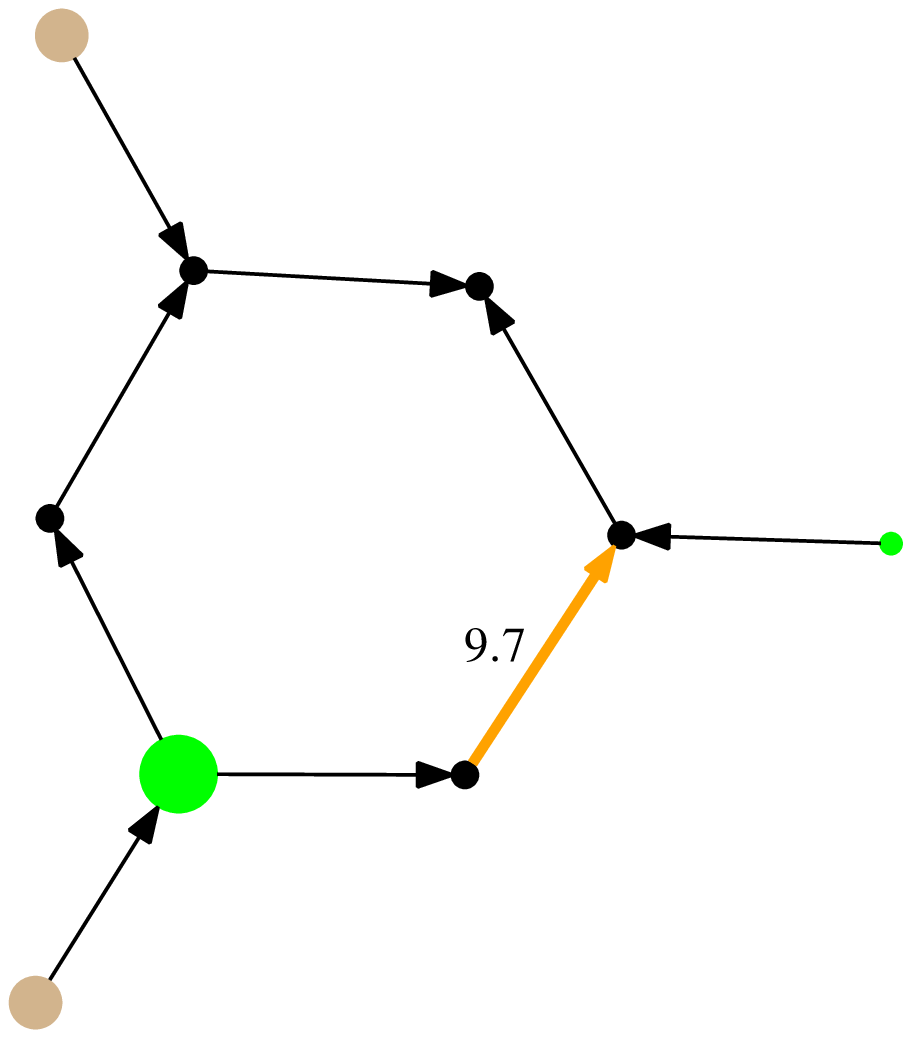}
\includegraphics[width=0.45\textwidth]{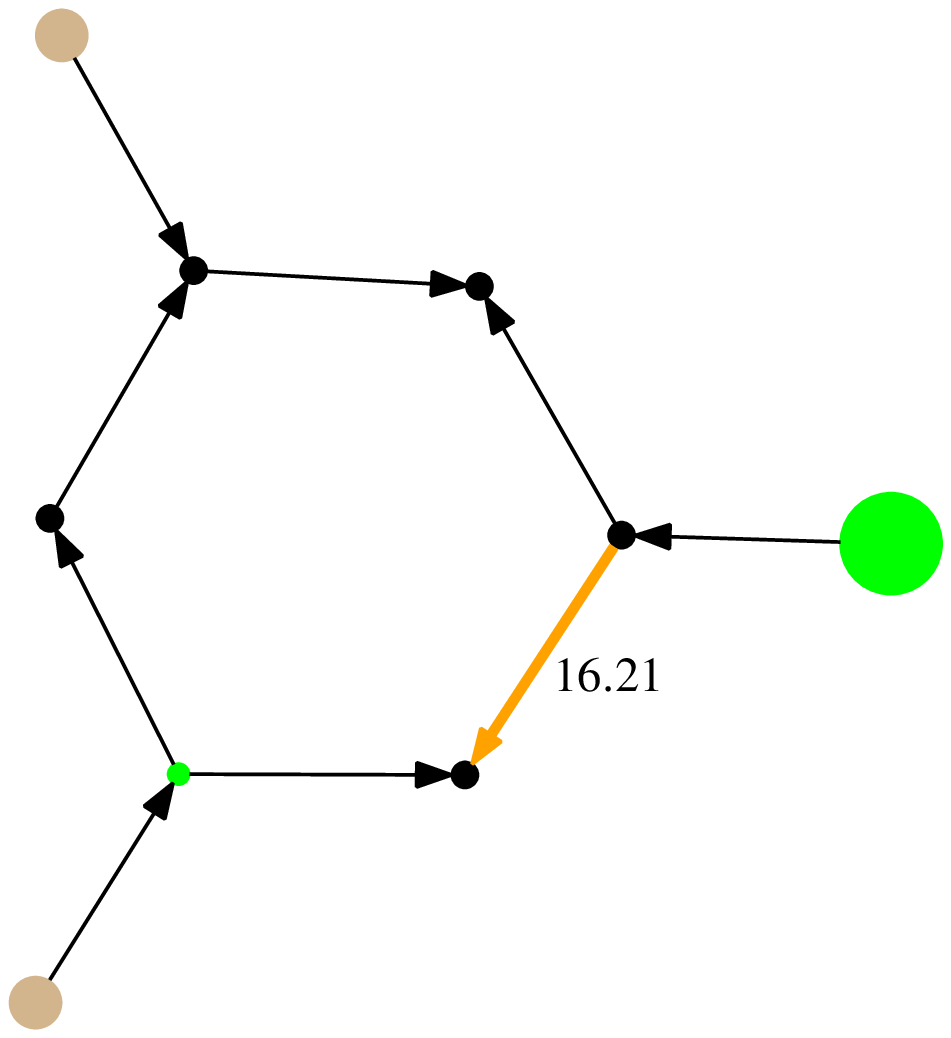} \caption{9-bus case,
25\% average penetration from two wind sources. With shifting winds, the flow
on the orange line changes direction with a large absolute
difference.\label{reversal}} \end{figure}

Here we consider the 9-bus case with two wind sources and 25\% average penetration and standard deviations set to 30\% of the average case and analyze the following two somewhat rare but still admissible wind configurations: (1) wind source (a) produces its average amount of power and source (b) three standard deviations \textit{below} average; (2) the reverse of the case (1). This results in a substantial reversal of flow on a particular line shown in Figure~\ref{reversal}. This example suggests that when evaluating and planning for grids with high-penetration of renewables one needs to be aware of potentially fast and significant structural rearrangements of power flows. Flow reversals and other qualitative changes of power flows, which are extremely rare in the grid of today, will become significantly much more frequent (typical) in the grid of tomorrow.

\subsection{Out-of-sample tests}
\label{subsec:out-of-sample}

We now study the performance of the control computed using \textit{nominal} CC-OPF when there are errors in the
underlying distribution of wind power. We consider two types of errors: (1)~the
true distribution is non-Gaussian but our Gaussian fit is ``close'' in an
appropriate sense, and (2)~the true distribution is Gaussian but with different
mean or standard deviation.  The experiments in this section use as data
set the BPA grid, which as noted before has 2209 buses and 2866 lines,
and collected wind data; altogether constituting a realistic test-case.

We first consider the non-Gaussian case, using the following probability
distributions, all with fatter tails than Gaussian: (1)~Laplace, (2)~logistic,
(3)~Weibull (three different shapes), (4)~t location-scale with 2.5 degrees of
freedom, (5)~Cauchy. For the Laplace and logistic distributions, we simply match
the mean and standard deviation. For the Weibull distribution, we consider shape
parameters $k=1.2,2,4$ and choose the scale parameter to match the standard
deviations. We then translate to match means. For the t distribution, we fix 2.5
degrees of freedom and then choose the location and scale to match mean and
standard deviation. For the Cauchy distribution, we set the location parameter
to the mean and choose the scale parameter so as to match the 95th percentiles.

We use our model and solve under the Gaussian assumption, seeking a solution
which results in no line violations for cases within two standard deviations of
the mean, i.e. a maximum of about 2.27\% chance of exceeding the limit. We then
perform Monte Carlo tests drawing from the above distributions to determine the
actual rates of violation. See Figure~\ref{oos1}. The worst-performer is the
highly-asymmetric (and perhaps unreasonable) Weibull with shape parameter 1.2,
which approximately doubles the desired maximum chance of overload. Somewhat
surprisingly, the fat-tailed logistic and Student's t distributions result in a
maximum chance of overload significantly less than desired, suggesting that our
model is too conservative in these cases.

\begin{figure}
\centering
\includegraphics[width=0.9\textwidth]{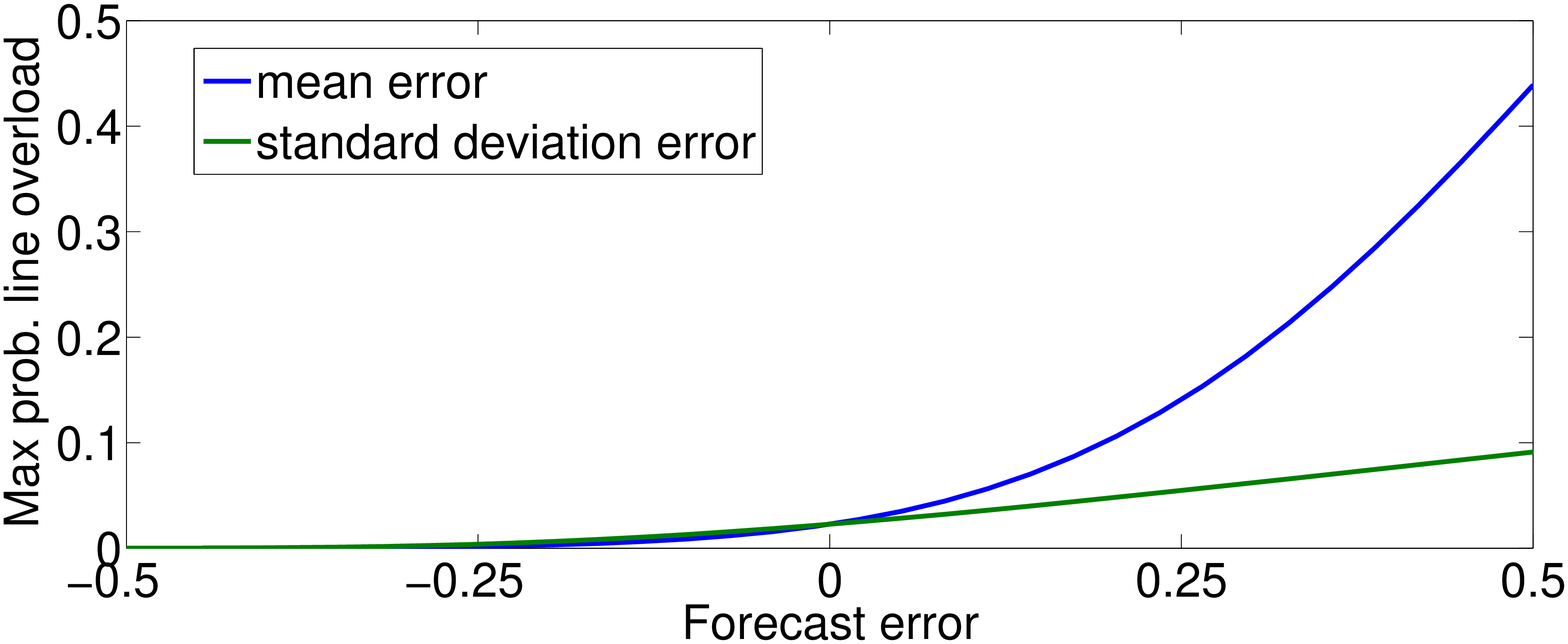}
\caption{BPA case solved with average penetration at 8\% and standard
deviations set to 30\% of mean. The maximum probability of line overload
desired is 2.27\%, which is achieved with 0 forecast error on the graph. Actual
wind power means are then scaled according to the x-axis and maximum
probability of line overload is recalculated (blue). The same is then done for
standard deviations (green).
\label{oos2}} \end{figure}

Next we consider the Gaussian case with errors. We solve with nominal values for
the mean and standard deviation of wind power. We then consider the rate of
violation after scaling all means and standard deviations (separately) . While
the solution is sensitive to errors in the mean forecast, the sensitivity is
well-behaved. With a desired safety level of $\epsilon=2.27\%$ for each line, an
error in the mean of 25\% results in a maximum 15\% chance of exceeding the
limit. The solution is quite robust to errors in the standard deviation
forecast, with a 25\% error resulting in less than 6\% chance of overload. See
Figure~\ref{oos2}.

\begin{figure}
\centering
\begin{tabular}{l c}
\toprule
Distribution & Max. prob. violation \\ \midrule
Normal & 0.0227 \\
Laplace & 0.0297 \\
logistic & 0.0132 \\
Weibull, $k=1.2$ & 0.0457 \\
Weibull, $k=2$ & 0.0355 \\
Weibull, $k=4$ & 0.0216 \\
t location-scale, $\nu=2.5$ & 0.0165 \\
Cauchy & 0.0276 \\
\bottomrule
\end{tabular}
\caption{ Maximum probability of overload for out-of-sample tests. These are a
result of Monte Carlo testing with 10,000 samples on the BPA case, solved under
the Gaussian assumption and desired maximum chance of overload at 2.27\%.
\label{oos1}} \end{figure}

\subsection{Scalability}
\label{subsec:scalability}

\begin{figure}[h]
\centering
\includegraphics[width=0.6\textwidth]{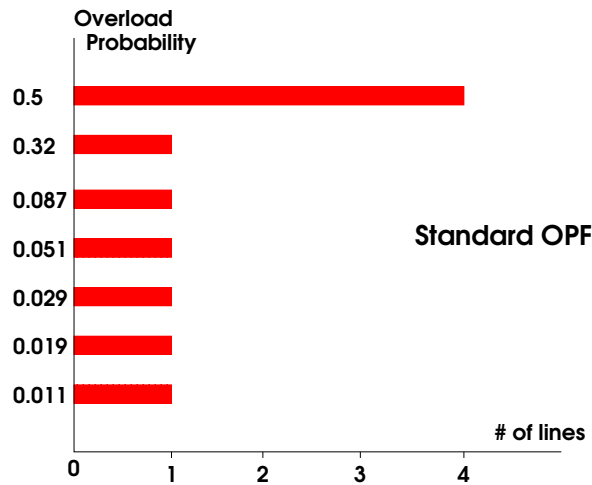}
\vskip 2pt
\caption{Figure shows number of lines that are overloaded with given probability values in simulation of 2746 bus Polish power grid using standard OPF with 20\% wind penetration distributed over $18$ wind farms. In particular, two lines are overloaded half of the time, and one line is overloaded one-third of the time, constituting a situation with unacceptable systemic risk. }\label{overloadsstandard}
\end{figure}

As an additional experiment illustrating scalability of the approach we studied the Polish national grid (obtained from MATPOWER as explained above) under simulated 20 \% renewable penetration spread over 18 wind farms, co-located with the 18 largest generators.  This co-location should lessen the risk associated with renewable fluctuation (which should be partially ``absorbed'' by the co-located generators). Figure \ref{overloadsstandard} studies the resulting risk exposure under standard OPF.  The chart shows the number of lines that attain several levels of overload probability.  The situation in the chart is unacceptable: it would lead to frequent tripping of at least four lines.

\begin{figure}[h]
\centering
\includegraphics[width=0.6\textwidth]{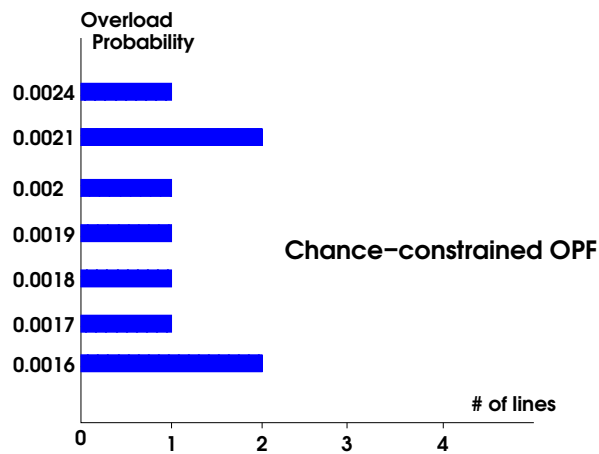}
\caption{Same as Figure \ref{overloadsstandard}, but under chance-constrained OPF. Notice that the largest overload probability is $200$ times smaller than in the case of standard OPF. Moreover, the cost increase is by less than one percent}\label{overloadscc}
\end{figure}

In contrast, Figure \ref{overloadscc} shows the performance attained by the chance-constrained OPF in the same setting as that of Figure \ref{overloadsstandard}.  Notice the drastic reduction in overload probabilities -- the system is stable.  Moreover, this is attained with a minor increase in cost (less than one percent) while the computational time is on the order of 10 seconds.


\section{Discussions}
\label{sec:discussions}

This manuscript suggests a new approach to incorporating uncertainty in the standard OPF setting routinely used in the power industry to set generation during a time window, or period (typically 15 min to one hour duration). When uncertainty associated with renewable generation is quantified in terms of the probability distribution of output during the next period, we incorporate it through chance constraints - probabilistic conditions which require that any line of the system will not be overloaded for all but a small fraction of time (at most one minute per hour, for example). Additionally, the modeling accounts for local frequency response of controllable generators to renewable changes. The key technical result of this manuscript is that the resulting optimization problem, CC-OPF, can be stated as a convex, deterministic optimization problem.  This result also relies on plausible assumptions regarding the exogenous uncertainty and linearity of the underlying power flow approximations/equations.  In fact, our CC-OPF is a convex (conic) optimization problem, which we solve very efficiently, even on realistic large-scale instances, using a sequential linear cutting plane algorithm.

This efficient CC-OPF algorithm becomes an instrument of our (numerical) experiments which were performed on a number of standard (and nonstandard) power grid data sets. Our experimental results are summarized as follows:
\begin{itemize}
\item We observe that CC-OPF delivers feasible results where standard OPF, run for the average forecast, would fail in the sense that many lines would be overloaded an unacceptably large portion of time.
\item Not only is CC-OPF safer than standard OPF, but it also results in \textit{cheaper} operation. This is demonstrated by  considering the optimal cost of CC-OPF under sufficiently high wind penetration solution (where standard OPF would fail) and the low penetration solution (corresponding to the highest possible penetration where standard OPF would not fail).
\item We discover that solutions produced by CC-OPF deviate significantly from what amounts to a the naive adjustment of the standard OPF obtained by correcting dispatch just at those generators which are close to overloaded lines.
\item  We test the level of wind penetration which can be tolerated without upgrading lines. This experiment illustrates that, at least for the model tested, the 30\% of wind penetration with rather strict probabilistic guarantees enforced by our CC-OPF may be feasible; but much lower wind penetration remains feasible under the standard approach.
\item We experiment with location of wind farms and discover strong sensitivity of the maximum level of penetration on the choice of location - optimal choice of wind farm location is critical for achieving the ambitious goal of high renewable penetration.
\item Analyzing fluctuations of line flows within CC-OPF solution admissible under high wind penetration, we discover that these fluctuations may be significant,  in particular resulting in reversal of power flows over some of the lines. This observation suggests that flow reversals and other qualitative changes of power flows, which are extremely rare in the grid of today, will become significantly much more frequent (typical) in the grid of tomorrow.
\item We studied an out-of-sample test consisted in applying CC-OPF (modeling exogenous fluctuations as Gaussian) to other distributions. Overall these tests suggest that with a proper calibration of the effective Gaussian distribution our CC-OPF delivers a rather good performance. One finds that the worst CC-OPF performance is observed for the most asymmetric distributions.
\item We also presented a computationally sound data-robust version of the CC-OPF where the parameters for the Gaussian distributions are assumed unknown, but lying in a window. This allows for parameter mis-estimation, for model error, and it is suggests a way to deal in a tractable way with non-Gaussian fluctuations.
\end{itemize}

The nature of the problem discussed in the manuscript -- the design of a new paradigm for computationally efficient generation re-dispatch that accounts for wind fluctuations -- inevitably required incorporation of a number of assumptions and approximations. In particular, we made simplifying assumptions about static forecasts and general validity of power flow linearization. We have also focused solely on failures associated with line congestion ignoring other possible difficulties, for example those associated with loss of synchronicity and voltage variations. However,  we would like to emphasize that all of these assumptions (admittedly natural for a first attack on the problem) also allow generalization within the approach just sketched:
\begin{itemize}
\item Accounting for time evolving forecast/loads/etc. Wind forecast, expressed in terms of the mean and standard deviation at the wind farm sites, changes on the scales comparable to duration of the generation re-dispatch interval. Loads may also change at these time scales. When the slowly evolving,  but still not constant, wind and load forecasts are available we may keep the quasi-static power flow description and incorporate this slow evolution in time into the chance constrained scheme. These changes will simply result in generalizing the conic formulation \eqref{obj}-\eqref{gen2} by splitting what used to be a single time interval into sub-intervals and allowing the regular generation to be re-dispatched and parallel coefficients to be adjusted more often.  Ramping rate constraints on the controllable generation may naturally be accounted in the temporal optimization scheme as well.
\item Accounting for nonlinearity in power flows. Evolution of the base case invalidates the linearization (DC-style) hypothesis. However, if variations around one base case becomes significant one may simply adjust the linearization procedure doing it not once (as in the case considered in the manuscript) but as often as needed.  Slow adjustment of the base case may also be included into the dynamical procedure mentioned one item above.  Additionally, some interesting new methodologies for handling nonlinearities have recently emerged, see \cite{12LL}.
\item Accounting for synchronization bounds. Loss of synchronicity and resulting disintegration of the grid is probably the most acute contingency which can possibly take place  in a power system. The prediction of those conditions under which the power grid will lose synchronicity is a difficult nonlinear and dynamic problem. However, as shown recently in \cite{12DCB}, one can formulate an accurate linear and static necessary condition for the loss of synchronicity.  A chance-constrained version of the linear synchronization conditions can be incorporated seamlessly in our CC-OPF framework.
\end{itemize}
Finally, we see many opportunities in utilizing the CC-OPF (possibly modified) as an elementary unit or an integral part of even more complex problems,  such as combined unit commitment (scheduling large power plants normally days, weeks or even months ahead) \cite{08WSL} with CC-OPF,  planning grid expansion \cite{Bent} while accounting for cost operation under uncertainty, or incorporating CC-OPF in mitigating emergency of possible cascades of outages \cite{05CTD,06NDKCL,10PTC,11PTC,11Bie_a,11Bie_b,12EH,12BHUZ}. In this context, it would be advantageous to speed up our already very efficient CC-OPF even further.  See, for example, \cite{danogarud}, \cite{danobook}. A different methodology, relying on distributed algorithms, can be found in  \cite{BP-Energy}.





\section{Acknowledgments}

We are thankful to S. Backhaus and R. Bent for their comments, and to K. Dvijotham for help with BPA data. The work at LANL was carried out under the auspices of the National Nuclear Security Administration of the U.S. Department of Energy at Los Alamos National Laboratory under Contract No. DE-AC52-06NA25396. MC and SH also acknowledge partial support of the DTRA Basic Research grant BRCALL08-Per3-D-2-0022.  DB and SH were partially supported by DOE award DE-SC0002676.




\bibliographystyle{siam}
\bibliography{Bib/robust,Bib/SmartGrid,Bib/CC-OPF,Bib/PR,Bib/sean}

\begin{thebibliography}{10}

\bibitem{12BH}
{\sc S.~Baghsorkhi and I.~Hiskens}, {\em Analysis tools for assessing the
  impact of wind power on weak grids}, SysCon,  (2012), pp.~1 --8.

\bibitem{nemben}
{\sc R.~Ben-Tal and A.~Nemirovski}, {\em Robust solutions of linear programming
  problems contaminated with uncertain data}, Mathematical Programming, 88
  (2000), pp.~411--424.

\bibitem{Bent}
{\sc R.~Bent, G.~Toole, and A.~Berscheid}, {\em Transmission network expansion
  planning with complex power flow models}, Power Systems, IEEE Transactions
  on, 27 (2012), pp.~904 --912.

\bibitem{00BV}
{\sc A.~Bergen and V.~Vittal}, {\em Power System Analysis}, Prentice Hall. Inc,
  2000.

\bibitem{12BHUZ}
{\sc A.~Bernstein, D.~Hay, M.~Uzunoglu, and G.~Zussman}, {\em Power grid
  vulnerability to geographically correlated failures - analysis and control
  implications}.
\newblock \url{http://arxiv.org/abs/1206.1099}, 2012.

\bibitem{bertsim}
{\sc D.~Bertsimas and M.~Sim}, {\em Price of robustness}, Operations Research,
  52 (2004), pp.~35--53.

\bibitem{danobook}
{\sc D.~Bienstock}, {\em Potential Function Methods for Approximately Solving
  Linear Programming Problems, Theory and Practice}, Kluwer, 2002.

\bibitem{11Bie_a}
{\sc D.~Bienstock}, {\em Adaptive online control of cascading blackouts}, Proc.
  2011 IEEE Power and Energy Society General Meeting,  (2011), pp.~1--8.

\bibitem{11Bie_b}
\leavevmode\vrule height 2pt depth -1.6pt width 23pt, {\em Optimal control of
  cascading power grid failures}, Proc 2011 CDC-ECC,  (2011), pp.~2166 --2173.

\bibitem{danogarud}
{\sc D.~Bienstock and G.~Iyengar}, {\em Faster approximation algorithms for
  covering and packing problems}, SIAM Journal on Computing, 35 (2006),
  pp.~825--854.

\bibitem{08BDL}
{\sc H.~Bludszuweit, J.~Dominguez-Navarro, and A.~Llombart}, {\em Statistical
  analysis of wind power forecast error}, Power Systems, IEEE Transactions on,
  23 (2008), pp.~983 --991.

\bibitem{BVbook}
{\sc S.~Boyd and L.~Vandenberghe}, {\em Convex Optimization}, Cambridge
  University Press, 2004.

\bibitem{BPA}
{\sc BPA-WIND}, {\em Wind generation \& total load in the bpa balancing
  authority}.
\newblock \url{http://transmission.bpa.gov/business/operations/wind/}, 2012.

\bibitem{CAISO-2007}
{\sc CAISO}, {\em Integration of renewable resources: Transmission and
  operating issues and recommendations for integrating renewable resources on
  the california iso-controlled grid}.
\newblock \url{http://www.caiso.com/1ca5/1ca5a7a026270.pdf}, 2007.

\bibitem{charnes}
{\sc A.~Charnes, W.~Cooper, and G.~Symonds}, {\em Cost horizons and certainty
  equivalents: an approach to stochastic programming of heating oil},
  Management Science, 4 (1958), pp.~235­--263.

\bibitem{05CTD}
{\sc J.~Chen, J.~Thorp, and I.~Dobson}, {\em Cascading dynamics and mitigation
  assessment in power system disturbances via a hidden failure model},
  International Journal of Electrical Power; Energy Systems, 27 (2005), pp.~318
  -- 326.

\bibitem{10CPS}
{\sc M.~Chertkov, F.~Pan, and M.~Stepanov}, {\em Predicting failures in power
  grids: The case of static overloads}, IEEE Transactions on Smart Grids, 2
  (2010), p.~150.

\bibitem{11CSPB}
{\sc M.~Chertkov, M.~Stepanov, F.~Pan, and R.~Baldick}, {\em Exact and
  efficient algorithm to discover extreme stochastic events in wind generation
  over transmission power grids}, CDC-ECC,  (2011), pp.~2174 --2180.

\bibitem{CIGRE09}
{\sc CIGRE}, {\em Technical brochure on grid integration of wind generation},
  International Conference on Large High Voltage Electric Systems,  (2009).

\bibitem{CPLEX}
{\sc CPLEX}, {\em {ILOG CPLEX Optimizer}}.
\newblock
  \url{http://www-01.ibm.com/software/integration/optimization/cplex-optimizer/},
  2012.

\bibitem{DENA}
{\sc DENA}, {\em Energy management planning for the integration of wind energy
  into the grid in germany, onshore and offshore by 2020},  (2005).

\bibitem{12DCB}
{\sc F.~{D{\"o}rfler}, M.~{Chertkov}, and F.~{Bullo}}, {\em {Synchronization in
  Complex Oscillator Networks and Smart Grids}}, Proceedings of National
  Academy of Sciences, 10.1073/pnas.1212134110 (2013).

\bibitem{20-2030}
{\sc EERE}, {\em 20\% of wind energy by 2030: Increasing wind energy's
  contribution to us electricity supply}, Department of Energy (DOE),
  DOE/GO-102008-2567 (2008).

\bibitem{12EH}
{\sc M.~Eppstein and P.~Hines}, {\em A "random chemistry" algorithm for
  identifying collections of multiple contingencies that initiate cascading
  failure}, Power Systems, IEEE Transactions on, 27 (2012), pp.~1698--1705.

\bibitem{ambiguous}
{\sc E.~Erdogan and G.~Iyengar}, {\em Ambiguous chance constrained problems and
  robust optimization}, Mathematical Programming, 107 (2007), pp.~37--61.

\bibitem{usc}
{\sc U.-C. P. S. O.~T. Force}, {\em Report on the august 14, 2003 blackout in
  the united states and canada: Causes and recommendations}.
\newblock \texttt{ https://reports.energy.gov}, (2004).

\bibitem{SOCP}
{\sc D.~Goldfarb and F.~Alizadeh}, {\em Second-order cone programming},
  Mathematical Programming, 95 (2003), pp.~3--51.

\bibitem{dongarud}
{\sc D.~Goldfarb and G.~Iyengar}, {\em Robust portfolio selection problems},
  Mathematical Programming, 28 (2001), pp.~1--38.

\bibitem{06Goz}
{\sc G.~Gonzalez and et~al.}, {\em Experience integrating and operating wind
  power in the peninsular spanish power system. point of view of the
  transmission system operator and a wind power producer}, CIGRE,  (2006).

\bibitem{GLS}
{\sc M.~Gr\"otschel, L.~Lov\'asz, and A.~Schrijver}, {\em Geometric Algorithms
  and Combinatorial Optimization}, Springer-Verlag, 1998.

\bibitem{gurobi}
{\sc GUROBI}, {\em {Optimizer}}.
\newblock \url{http://www.gurobi.com/}, 2012.

\bibitem{nrelweibull}
{\sc B.-M. Hodge and M.~Milligan}, {\em Wind power forecasting error
  distributionsover multiple timescales}, Detroit Power Engineering Society
  Meeting,  (2011), pp.~1--8.

\bibitem{Huneault1991}
{\sc M.~Huneault}, {\em {A survey of the optimal power flow literature}}, Power
  Systems, IEEE Transactions,  (1991).

\bibitem{BP-Energy}
{\sc M.~Kraning, E.~Chu, J.~Lavaei, and S.~Boyd}, {\em Message passing for
  dynamic network energy management}.
\newblock \url{http://www.stanford.edu/~boyd/papers/decen_dyn_opt.html}, 2012.

\bibitem{94Kun}
{\sc P.~Kundur}, {\em Power System Stability and Control}, McGraw-Hill, New
  York, NY, USA, 1994.

\bibitem{12LL}
{\sc J.~Lavaei and S.~Low}, {\em Zero duality gap in optimal power flow
  problem}, IEEE Transactions on Power Systems, 27 (2012), pp.~92--107.

\bibitem{Makarov-wind}
{\sc Y.~Makarov, C.~Loutan, J.~Ma, and P.~de~Mello}, {\em Operational impacts
  of wind generation on california power systems}, Power Systems, IEEE
  Transactions on, 24 (2009), pp.~1039 --1050.

\bibitem{millerwagner}
{\sc L.~Miller and H.~Wagner}, {\em Chance-constrained programming with joint
  constraints}, Operations Research, 13 (1965), pp.~930--945.

\bibitem{mosek}
{\sc MOSEK}, {\em {Optimizer}}.
\newblock \url{http://www.mosek.com/}, 2012.

\bibitem{06NDKCL}
{\sc D.~Nedic, I.~Dobson, D.~Kirschen, B.~Carreras, and V.~E. Lynch}, {\em
  Criticality in a cascading failure blackout model}, International Journal of
  Electrical Power; Energy Systems, 28 (2006), pp.~627 -- 633.

\bibitem{06NS}
{\sc A.~Nemirovski and A.~Shapiro}, {\em Convex approximations of chance
  constrained programs}, SIAM Journal on Optimization, 17 (2006), pp.~969--996.

\bibitem{Ozturk2004}
{\sc U.~Ozturk and M.~Mazumdar}, {\em {A solution to the stochastic unit
  commitment problem using chance constrained programming}}, Power Systems,
  IEEE, 19 (2004), pp.~1589--1598.

\bibitem{11PTC}
{\sc R.~{Pfitzner}, K.~{Turitsyn}, and M.~{Chertkov}}, {\em {Controlled
  Tripping of Overheated Lines Mitigates Power Outages}}, arxiv:1104.4558,
  (2011).

\bibitem{10PTC}
{\sc R.~Pfitzner, K.~Turitsyn, and M.~Chertkov}, {\em Statistical
  classification of cascading failures in power grids}, in Power and Energy
  Society General Meeting, 2011 IEEE, july 2011, pp.~1 --8.

\bibitem{chance}
{\sc A.~Pr\'{e}kopa}, {\em On probabilistic constrained programming}, Proc.
  Princeton Symposium on Math. Programming,  (1970), pp.~113--138.

\bibitem{12VMLA}
{\sc M.~Vrakopoulou, K.~Margellos, J.~Lygeros, and G.~Andersson}, {\em A
  probabilistic framework for security constrained reserve scheduling of
  networks with wind power generation}, in Energy Conference and Exhibition
  (ENERGYCON), 2012 IEEE International, sept. 2012, pp.~452 --457.

\bibitem{meyn-markets-volatility}
{\sc G.~Wang, M.~Negrete-Pincetic, A.~Kowli, E.~Shafieepoorfard, S.~Meyn, and
  U.~Shanbhag}, {\em Dynamic competitive equilibria in electricity markets}, in
  Control and Optimization Theory for Electric Smart Grids, A.~Chakrabortty and
  M.~Illic, eds., Springer, 2011.

\bibitem{08WSL}
{\sc J.~Wang, M.~Shahidehpour, and Z.~Li}, {\em Security-constrained unit
  commitment with volatile wind power generation}, Power Systems, IEEE
  Transactions on, 23 (2008), pp.~1319 --1327.

\bibitem{NYT2008}
{\sc M.~Wyld}, {\em Wind energy bumps into power grid limits}, New York Times,
  (2008/08/27).

\bibitem{Zhang2011}
{\sc H.~Zhang}, {\em {Chance Constrained Programming for Optimal Power Flow
  Under Uncertainty}}, Power Systems, IEEE Transactions on, 26 (2011),
  pp.~2417--2424.

\bibitem{11ZSHY}
{\sc S.~Zhang, Y.~Song, Z.~Hu, and L.~Yao}, {\em Robust optimization method
  based on scenario analysis for unit commitment considering wind
  uncertainties}, in Power and Energy Society General Meeting, 2011 IEEE, july
  2011, pp.~1 --7.

\bibitem{MATPOWER}
{\sc R.~Zimmerman, C.~Murillo-Sanchez, and R.~Thomas}, {\em Matpower:
  Steady-state operations, planning, and analysis tools for power systems
  research and education}, Power Systems, IEEE Transactions on, 26 (2011),
  pp.~12 --19.

\end{thebibliography}

\end{document}